\documentclass[reqno, 12pt]{amsart}
\def\C{{\mathbb C}}
\def\Cp{{\mathbb C}_p}

\def\Qp{{\mathbb Q}_p}
\def\Z{{\mathbb Z}}
\def\Zp{{\mathbb Z}_p}
\def\R{{\mathbb R}}

\def\Fp{{{\mathbb F}_p}}
\def\F{{\mathbb F}}
\def\Qp {{{\mathbb Q}_p}}
\def\Fq{{{\mathbb F}_q}}
\def\Fr{{{\mathbb F}_r}}

\def\dis{\displaystyle}

\newtheorem{lemma}{Lemma}
\newtheorem{cor}{Corollary}
\newtheorem{prop}{Proposition}

\newtheorem{theorem}{Theorem}

\theoremstyle{definition}
\newtheorem{defn}{Definition}   

\newtheorem{question}{Question}

\theoremstyle{remark}
\newtheorem{rem}{Remark}        
\newtheorem{rems}{Remarks}      
\newtheorem{example}{Example}

\begin{document}
\title[Applications of non-Archimedean integration to the $L$-series of
$\tau$-sheaves]
{Applications of non-Archimedean integration to the $L$-series of 
$\tau$-sheaves}
\author{David Goss}
\thanks{This work is dedicated to the memory of Arnold Ross; a wonderful
colleague who had an enormous positive influence on the mathematics community}
\address{Department of Mathematics\\ The Ohio State University\\ 231 W.
$18^{\text{th}}$ Ave. \\ Columbus, Ohio 43210}
\email{goss@math.ohio-state.edu}
\date{April 30, 2004}

\begin{abstract}
Let $\underline{\mathcal F}$ be a $\tau$-sheaf. Building on previous
work of Drinfeld, Anderson, Taguchi, and Wan, B\"ockle and Pink \cite{bp1}
develop a  cohomology theory for $\underline{\mathcal F}$. In \cite{boc1}
B\"ockle uses this theory to establish the analytic continuation of
the $L$-series 
associated to $\underline{\mathcal F}$ (which is a characteristic $p$ valued 
``Dirichlet series'') {\em and} the logarithmic
growth of the degrees of its special polynomials. In this
paper we shall show that this logarithmic growth is all that is needed to
analytically continue the original $L$-series as well as {\em all}
associated partial $L$-series. Moreover, we show that the degrees of
the special polynomials attached to the partial $L$-series also grow
logarithmically. Our tools are B\"ockle's original results, non-Archimedean
integration, and the very strong estimates of Y.\ Amice \cite{am1}. Along the
way, we define certain natural modules associated with non-Archimedean
measures (in the characteristic $0$ case as well as in characteristic $p$).

\end{abstract}

\maketitle

\section{Introduction}\label{intro}
In his original work \cite{r1} on his zeta function, Riemann established
that the density of zeroes up to level $T$ in the critical strip is
approximately $\frac{1}{2\pi}\log(\frac{T}{2\pi})$. Since then similar
results have been established for general $L$-series.

In the arithmetic of function fields over finite fields, logarithmic
growth manifests itself for characteristic $p$ valued zeta functions in 
terms of the degrees of their associated ``special polynomials''
(see Subsection \ref{entire}). This was first noted by the author in
the explicit measure calculations of Dinesh Thakur \cite{th1} 
for $\Fr[T]$.

More generally,
let $k$ be an arbitrary global function field with full field of constants
$\Fr$ and let $\infty$ be a fixed place. Set 
$A$ to be the Dedekind domain of elements of $k$ which are integral
at all places outside $\infty$. Let $k_1$ be a finite extension of $k$ and let
$\phi$ be a Drinfeld module over $k_1$. As in \S 8.6 of \cite{go4},
one can define the $L$-series $L(s)$ of $\phi$ which is a ``Dirichlet series''
$\sum_I c_I I^{-s}$, $I$ an ideal of $A$, in finite characteristic. Using
elementary estimates, (Lemma 8.8.1 of \cite{go4}) it was shown
when $\phi$ has rank $1$ (or complex multiplication etc.)
that $L(s)$ has an analytic continuation to an entire function.
Moreover analytic
continuations can then also be established for the interpolations
associated to $L(s)$ at finite primes in the sense of Subsection
\ref{entire}.

The estimates of \cite{go4} allow one to establish the existence of
special polynomials in the general rank $1$ case but give poor estimates on
the degrees of these polynomials. As mentioned, the explicit calculation of
Thakur in \cite{th1}, as well as the calculations of Newton polygons
by Wan, Diaz-Vargas and Sheats (\cite{wa1}, \cite{dv1} and \cite{sh1}), show
that these degrees in fact grow logarithmically and that this logarithmic
growth reflect rationality (in terms of the complete field $k_\infty$)
of the zeroes of such function. 
Moreover, this logarithmic growth, when combined with the deep a-priori
estimates of Amice \cite{am1}, actually provides the analytic continuation
of these $L$-series at {\em all} places of $k$.
It thus became reasonable to also expect
such logarithmic growth for the degrees of special polynomials
associated to $L$-series of general rank Drinfeld
modules, $t$-modules etc.  And, indeed, such a basic result
was recently established by G.\ B\"ockle in \cite{boc1} as a stepping
stone in his analytic continuation of such $L$-series.

The approach to characteristic $p$ $L$-series in \cite{boc1} is via
cohomology. Due to the labors of Drinfeld, G.\ Anderson, Y.\ Taguchi,
D.\ Wan, R.\ Pink and B\"ockle,
 the basic construction of characteristic $p$ arithmetic
has evolved from Drinfeld modules to ``$\tau$-sheaves'' which are
simply coherent sheaves over the product of a base scheme $X$ with
$\text{Spec}(A)$ equipped with a Frobenius-linear morphism $\tau$
(see Definition \ref{tausheaf}).
In \cite{bp1}, Pink and B\"ockle show how to embed the $\tau$-sheaves
into a category of ``crystals'' which possesses a good cohomology theory
and Lefschetz fixed point theorem. It is this cohomology theory that
B\"ockle uses to establish the logarithmic growth  of the
degrees of the special polynomials (and the analytic
continuation of these functions at $\infty$ and {\em all}
interpolations at finite primes) in very great generality.

Let $L(s)$ be the $L$-series of a $\tau$-sheaf of the type
shown to be entire in \cite{boc1} (see, e.g., Theorem \ref{cohenmc1} below).
 In this paper we show how the logarithmic growth of the degrees of
$L(x-j)$ is enough to establish the analytic continuation and logarithmic
growth of {\it any} {\it partial $L$-series} (see Definition \ref{partialL})
associated to $L(s)$. We say that such a Dirichlet series 
is ``in the motivic
class $\mathcal M$'' as it is our expectation that the only way to
provide non-trivial examples is precisely via $\tau$-sheaves.
The idea behind the proof is again to express $L(s)$
as a non-Archimedean integral which a priori is valid in the ``half-plane''
where $L(s)$ has absolute convergence (in the sense
of Remark \ref{absconv}). Then one uses the logarithmic growth of the
degrees of
$L(x,-j)$ to show that the measures so obtained blow up very slowly.
Combined with Amice's estimates, the analytic continuation follows.
Then playing off Amice's estimates against certain a priori estimates on
coefficients 
allows one to also obtain the logarithmic growth for the degrees of the 
special polynomials associated to the partial $L$-series. 

A very interesting feature of the proof of the logarithmic growth of
the degrees of
special polynomials associated to partial $L$-series
 is the way that the theories at the place $\infty$
and {\em all} finite places intertwine. Indeed, to establish the logarithmic
growth of a partial $L$-series defined modulo a given place $w$ we
need crucially to use the $w$-adic theory associated to $L(s)$. 

Along the way, we elucidate some of the formalism associated with
non-Archimedean integration both in finite characteristic and
characteristic $0$. In particular we show how the convolution 
product of measures comes equipped with certain canonical associated
modules. In finite characteristic, these modules give a concrete
realization of the space of measures as ``differential
operators'' which was previously only known abstractly (e.g., \cite{go8}).

Our results, along with those of B\"ockle, Pink, Taguchi, and Wan,
make it very reasonable to hope that a deeper theory of the zeroes
will eventually be found. Indeed, the results of \cite{wa1}, \cite{dv1},
\cite{sh1} give far more information in certain special cases than is
obtainable from the estimates given here.
As of now this theory would seem to
involve first a deeper understanding of the relationship between
the characteristic $p$ $L$-series and modular forms associated to
Drinfeld modules as established in \cite{boc2} (and presented in
\cite{go7}). Indeed, B\"ockle in \cite{boc2} associates
a $\tau$-sheaf to a cusp-form via Hecke operators; thus cusp-forms 
{\em also} give rise to characteristic $p$ valued $L$-series.

The reader may wonder why one could not approach our results by simply
using twists of the $L$-series by abelian characters and then solving
for the partial $L$-series. However, there are simply not enough
characters with values in finite characteristic for this to work in
general.

In this paper we have worked with completely arbitrary $A$. All the results
go through in this case, but there are a number of associated technicalities
that must be dealt with. These technicalities involve making
sense of ``$I^s$'' when $I$ is not generated by a ``monic
element'' (in
the sense of Subsection \ref{domaininfty} which generalizes the usual
notion of monic polynomial). Of course
when $A=\Fr[T]$, all $I$ obviously {\em are} so generated, and thus
the technicalities vanish.
Therefore the reader is well advised to first read this
paper with $A=\Fr[T]$.

It is my pleasure to thank Gebhard B\"ockle for his help in understanding
the results in \cite{bp1}, \cite{boc1}, and \cite{boc2} and for his
comments on early versions of this work. These comments greatly helped clarify
the proof of our main result.
Similarly,  I am also indebted to Keith Conrad and the
referee for helpful comments.
It is finally 
my pleasure to thank Zifeng Yang who pointed out that my original
proof of the logarithmic growth (of the degrees of the special polynomials)
at the infinite place would also work at the finite places. Indeed,
because the degree of a principal divisor on a complete curve must equal $0$,
the order of zero of an element of $A$ at a prime of $A$ is obviously
controlled by the order of its pole at $\infty$. 

\section{Review of non-Archimedean integration}

\subsection{General theory}\label{gentheory}

In this section
$K$ will be a non-Archimedean local field, of any characteristic,
with maximal compact subring $R_K$ and associated maximal ideal $M_K$. 
Thus $\F_K:= R_K/M_K$ is
a finite field and we denote its order by $q_{_K}$. Let $\vert ?\vert=\vert ?\vert_K$ be the absolute value on $K$ 
defined by $\vert \pi\vert=q_{_K}^{-1}$, where
$M_K=(\pi)$ (so that
$R_K=\{x\in K\mid \vert x\vert \leq 1\}$ and $M_K=\{x\in K\mid
\vert x\vert <1\}$). We let $v(?)=v_K(?)$ be the 
additive valuation associated to
$\vert ?\vert$ with $v(\pi)=1$.
Let $\bar{K}$ be a fixed 
algebraic closure of $K$ equipped with the canonical extension of 
$\vert ?\vert$ and $v(?)$. Finally let $K^{\rm sep}\subseteq \bar{K}$ be the
separable closure.

Let $L\subset \bar{K}$ be a finite extension of $K$ 
with integers $R_L$ (so that
$L$ is still a local field).

\begin{defn} \label{measure}
An {\it $R_L$-valued measure}
on $R_K$ is a finitely additive $R_L$-valued function
on the compact open subsets of $R_K$.\end{defn}

\noindent
More generally, one defines an {\it $L$-valued measure on $R_K$} to be a finitely
additive $L$-valued function $\mu$ on the compact open subsets of $R_K$ with
bounded image in $L$.
One sees immediately that the $\mu$ is an $L$-valued measure if and
only if there exists $a\neq 0\in L$ such that $a\mu$ is an
$R_L$-valued measure. We will denote
the space of $R_L$-valued measures on $R_K$ by ${\mathcal M}(R_K,R_L)$ and
the space of $L$-valued measures by ${\mathcal M}(R_K,L)$; so
${\mathcal M}(R_K,L)=L\otimes {\mathcal M}(R_K,R_L)$.

\begin{rem}\label{distribution}
Arbitrary (i.e., possibly unbounded)
$L$-valued finitely additive functions on the compact opens of
$R_K$ are called {\it $L$-valued distributions}.
\end{rem}

Let $f\colon R_K\to L$ be a continuous function and let 
$\mu\in {\mathcal M}(R_K,L)$. One defines
Riemann sums associated to $f$ and $\mu$ in the obvious manner. As
$R_K$ is compact, $f$ is also
uniformly continuous. Therefore it is easy to see that the Riemann
sums converge to an element of $L$ which is naturally denoted
$\displaystyle \int_{R_K} f(x)\, d\mu(x)$.

Let $E$ be a vector space over $L$. A map $\Vert~\Vert\colon E\to \R$ is
a {\it norm} if and only if
\begin{enumerate}
\item $\Vert x\Vert =0 \Leftrightarrow x=0\in E$,
\item $\Vert x+y\Vert\leq \max\{\Vert x\Vert, \Vert y\Vert\},$
\item $\Vert a x\Vert=\vert a\vert \Vert x\Vert$ for $a\in L$ and $x\in E$.
\end{enumerate}
The norm $\Vert ?\Vert$ induces an ultrametric $\rho$ on $E$
by $\rho(x,y):=\Vert x-y\Vert$. 
\begin{defn}\label{banachspace}
A {\it Banach space over $L$} is a complete normed $L$-vector space.
\end{defn}

Let $E$ be an $L$-Banach space. We say that $E$ is {\it separable}
if and only if $E$ contains a dense countable subset. From now on we will
only consider separable $L$-Banach spaces.

\begin{defn}\label{orthobasis}
Let $E$ be a Banach space and $\{e_i\}_{i=0}^\infty$ be a countable
subset of $E$. We say that $\{e_i\}$ is an {\it orthonormal} basis 
(or {\it Banach basis}) for $E$ if and only if
\begin{enumerate}
\item every $x\in E$ can be written uniquely as a convergent sum 
$x=\sum_{i=0}^\infty c_ie_i$ for $\{c_i\}\subset K$, $c_i\to 0$ as
$i\to \infty$, 
\item $\Vert x\Vert =\sup_i \{c_i\}$.
\end{enumerate}
\end{defn}

\begin{example}\label{contfunctions}
Let ${\mathcal C}(R_K,L)$ be the $L$-vector space of continuous $L$-valued functions
on $R_K$ equipped with the sup norm $\Vert f\Vert$ for continuous
$f$; it is easy to see that ${\mathcal C}(R_K,L)$ is an $L$-Banach space. 
A theorem due to Kaplansky (\cite{ka1} or Th.\ 5.28 in \cite{vr1}) assures
us that the polynomial functions are dense in ${\mathcal C}(R_K,L)$. It follows
readily that ${\mathcal C}(R_K,L)$ is separable.
\end{example}

The existence of orthonormal bases $\{Q_n(x)\}$, where $Q_n(x)$ is a polynomial
of degree $n$, for
${\mathcal C}(R_K,L)$ (where $L$ is local as above)
will be of critical importance to us. In \cite{am1}
(see also \cite{ya1}) Y.\ Amice constructs such a basis by first using Newton 
interpolation involving certain ``very well distributed'' sequences of 
elements of $R_K$ to construct polynomials $\{p_n(x)\}$ of degree
$n$. Then the Banach basis $\{Q_n(x)\}$ 
for ${\mathcal C}(R_K,L)$ is defined by $Q_n(x):=p_n(x)/s_n$ for all
$n$ where
$v_K(s_n)=\sum_{i=1}^\infty [n/q_K^i]$. 
In the case where
$R_K=\Fr[[T]]$, K.\ Conrad in \cite{co1} shows
how to use an orthonormal basis of the Banach space ${\mathcal L}_{\Fr}(R_K,L)$
of all {\it $\Fr$-linear} continuous functions (obviously a closed
subspace of ${\mathcal C}(R_K,L)$)
to construct a polynomial basis for the space of all continuous functions
via the ``digit principle.''
We shall have more to say about this later in Subsection \ref{hyperder}.

Let $\mu$ be a measure on $R_K$. Set
\begin{equation}\label{meascoeff}
\left\{b_n:=\int_{R_K} Q_n(x)\, d\mu(x)\,,\quad n=0,2\ldots\right\}\,.\end{equation}
We call $\{b_n\}$ the {\it measure coefficients associated to the basis
$\{Q_n(x)\}$}.
The boundedness of $\mu$ immediately implies that $\{b_n\}\subset L$ is also
bounded. Moreover, let $f(x)=\sum_{n=0}^\infty a_n Q_n(x)$ be a
continuous function, where the
{\it expansion coefficients} $\{a_n\}$ lie in $L$ and
$a_n\to 0$ as $n\to \infty$. Continuity implies that
\begin{equation}\label{intdef}
\int_{R_K} f(x)\, d\mu(x)=\sum_{n=0}^\infty a_nb_n\,.
\end{equation}
Note that the locally constant functions are also continuous; thus
{\it any} bounded sequence $\{b_n\}\subset L$ gives rise to
a bounded measure by Equation \ref{intdef}. Consequently a given choice of
orthonormal basis for ${\mathcal C}(R_K,L)$ immediately gives a corresponding
isomorphism of the space of measures with the space of bounded sequences
with coefficients in $L$.

\begin{defn}\label{dirac}
Let $\alpha\in R_K$ and $f\in {\mathcal C}(R_K,L)$. We define
the {\it Dirac measure} $\delta_\alpha$ associated to $\alpha$ by
$$\int_{R_K} f(x)\, d\delta_\alpha (x):=f(\alpha)\,.$$
\end{defn}
\noindent
Using the ideas just presented, one sees readily that the Dirac measure
is indeed a bounded measure in the sense of Definition \ref{measure}.
The Dirac measures
provide the basic building blocks for the constructions given in this
paper.

The space of measures is also an $L$-algebra via the convolution product
in the standard manner which we recall in the next definition.

\begin{defn} \label{convolution}
Let $\mu$ and $\nu$ be elements of ${\mathcal M}(R_K,L)$. Let
$f\colon R_K\to L$ be continuous. We define the
{\it convolution} $\mu\ast \nu\in {\mathcal M}
(R_K,L)$ by
\begin{equation}\label{convolution2}
\int_{R_K}f(u)\, d\mu\ast \nu(u):=\int_{R_K}\int_{R_K} f(x+y)\, d\mu(y)
d\nu(x)\,.\end{equation}
\end{defn}
\noindent
It is easy to see that Equation \ref{convolution2} does indeed define
a new measure and makes ${\mathcal M}(R_K,L)$ into a commutative
$L$-algebra.

Let $\alpha$ and $\beta$ be two elements of $R_K$. By definition one
has
\begin{equation}\label{convoldirac}
\delta_\alpha \ast \delta_\beta=\delta_{\alpha+\beta}\,.
\end{equation}
\subsection{The characterization of locally analytic functions}\label{locan}
In this subsection we will review the basic results of Y.\ Amice
\cite{am1} (see also \cite{ya1}) that permit
us to characterize those $f\in {\mathcal C}(R_K,L)$ which are locally
analytic.

Let $0\neq \rho\in R_K$ with $t=\vert \rho \vert$. Let $\alpha$ be another
element of $R_K$. The closed ball $B_{\alpha,t}$ around $\alpha$ of 
radius $t$ is defined, as usual, by
$$B_{\alpha,t}:=\{x\in R_K\mid \vert x-\alpha\vert \leq t\}\,.$$

\begin{defn}\label{analB}
A continuous function $f\colon B_{\alpha,t}\to L$ is {\it analytic on 
$B_{\alpha,t}$} if and only if $f$ may be expressed 
as 
\begin{equation}\label{analBpower}
f(x)=\sum_{n=0}^\infty b_n \left(\frac{x-\alpha}{\rho}\right)^n\,,
\end{equation}
where $\{b_n\}\subset L$ and $b_n\to 0$ as $n\to \infty$.
\end{defn}

\noindent
The norm $\Vert f\Vert_B$ of $f$ on $B:=B_{\alpha,t}$ is defined by
\begin{equation}\label{norm1}
\Vert f\Vert_B:=\sup_{n=0}^\infty \{\vert b_n\vert\}\,.\end{equation}
Clearly the set of functions analytic on $B$ forms an algebra which is
topologically isomorphic to the Tate algebra 
$L\langle \langle u
\rangle \rangle$ of power series $\sum_{i=0}^\infty c_iu^i$ converging on
the closed unit disc (i.e., 
$c_i\to 0$ as $i\to \infty$). 
Let $\bar{R}_K\subset \bar{K}$ be the ring of integers. Standard
results on Tate algebras then imply the basic result
\begin{equation}\label{tatenorm}
\Vert f \Vert_B=\sup_{\lambda\in \alpha+\rho \bar{R}_K}\{\vert f(\lambda)\vert\}\,.
\end{equation}

\begin{defn}\label{locanfunc}
A continuous function $f\in {\mathcal C}(R_K,L)$ is said to be
{\it locally analytic} if for each $\alpha\in R_K$ 
there exists $t_\alpha > 0$ such that $f$ is analytic on $B_{\alpha,t_\alpha}$.
\end{defn}

Now let $\pi\in R_K$ be a uniformizing parameter. 
\begin{defn}\label{orderh}
We say that a locally analytic function $f$ has {\it order $h$}, where
$h$ is a non-negative integer, if we can choose
$t_\alpha \geq \vert \pi \vert^h$ for all $\alpha\in R_K$. 
\end{defn} 

\noindent Definition \ref{orderh} is equivalent to requiring $f$ to be
analytic on each coset of $M_K^h$.
Note that by compactness one can find a finite number of $\alpha$ with
$\{B_{\alpha,t_\alpha}\}$ covering $R_K$. Thus every locally analytic function has
order $h$ for some non-negative $h$. 

We denote the space of locally analytic $L$-valued functions on $R_K$ by 
${\mathcal LA}(R_K,L)$
and those of order $h$ by ${\mathcal LA}_h(R_K,L)$. Clearly ${\mathcal LA}(R_K,L)
=\bigcup {\mathcal LA}_h(R_K,L)$.

\begin{defn}\label{locanalbanach}
Let $R_K=\bigcup_{j=0}^m B_j$ where the balls 
$B_j:= B_{\alpha_j, \vert \pi\vert^h}$ 
are mutually disjoint. Let $f\in {\mathcal LA}_h(R_K,L)$. Then we set
\begin{equation}\label{norm2}
\Vert f\Vert_h:=\max_j \{\Vert f\Vert_{B_j}\}\,.
\end{equation}
\end{defn}

\noindent One checks easily that Definition \ref{locanalbanach} makes
the space ${\mathcal LA}_h(R_K,L)$ a Banach space. One also readily sees that
a sequence of functions $\{f_i\}$ converging to a function $f$ in
${\mathcal LA}_h(R_K,L)$ will also converge to $f$ in ${\mathcal LA}_{h^\prime}
(R_K,L)$ for any $h^\prime \geq h$. By Equation \ref{tatenorm}, $\{f_i\}$ also
converges to $f$ in the sup norm on continuous functions with domain $R_K$.

Let $R_K=\bigcup B_j$ as in Definition \ref{locanalbanach}. Let $\{\chi_{j,n}(x)\}$,
$j=0,\cdots, m$ and $n\geq 0$, be the set of locally analytic functions defined
by 
\begin{equation}\label{chibasis}
\chi_{j,n}(x):=\begin{cases} \left(\frac{x-\alpha_j}{\pi^h}\right)^n &
\text{for}~
  x\in B_j\\
                0& \text{otherwise}\,.\end{cases}
\end{equation}
It is very simple to see that $\{\chi_{j,n}(x)\}$ is an orthonormal basis
for ${\mathcal LA}_h(R_K,L)$; thus ${\mathcal LA}_h(R_K,L)$ is also
a separable Banach space.

Let $\{Q_n(x)=p_n(x)/s_n\}$ be the orthonormal
basis for ${\mathcal C}(R_K,L)$ constructed by Amice as mentioned
above, and let
$f\in {\mathcal C}(R_K,L)$ be expressed as
$f(x)=\sum_{n=0}^\infty a_n Q_n(x)$ where $\{a_n\}\subset L$ and
$a_n\to 0$ as $n\to \infty$. Put
$$\gamma=\gamma_f:=\liminf_n \{v(a_n)/n\}\,.$$
We then have the following results \cite{am1} (see also \cite{ya1}).

\begin{theorem}\label{amice1}
{\rm 1.}\ The set $\{p_n(x)\}$ forms an orthonormal basis for the Tate
algebra of locally analytic functions of order $0$ (i.e., functions
analytic on the closed unit disc). Moreover, for $h\geq 1$ the collection
$\{p_n(x)/s_{n,h}\}$, $\{s_{n,h}\}\subset L$, forms an orthonormal
basis for ${\mathcal LA}_h(R_K,L)$ if and only if 
\begin{equation}\label{valsnh}
v_K(s_{n,h})=\sum_{i=1}^h [n/q_K^i]\,,\end{equation}
(where $[?]$ is the standard greatest integer function).\\
{\rm 2}.\ The function $f$ is locally analytic of order $h$ if and only if
$v(a_n)-\sum_{i=h+1}^\infty \left[\frac{n}{q^i}\right] $ tends to
$\infty$ as $n\to \infty$.\\
{\rm 3}.\ The function $f$ is locally analytic if and only if $\gamma >0$.
In this case, set $$l:=\max\{0,1+[-\log(\gamma (q-1))/\log q]\}\,.$$ Then
$f$ is locally analytic of order $h\geq l$.
\end{theorem}

\noindent
Note that Part 1 of Theorem \ref{amice1} can easily be restated in
terms of $\{Q_n(x)\}$. Using this reformulation,
in \cite{ya1}, Z.\ Yang shows that Theorem \ref{amice1} remains true for
{\it any} polynomial orthonormal
basis $\{h_n(x)\}$ for ${\mathcal C}(R_K,L)$ with
$\deg h_n(x)=n$ all $n$. In particular, it holds true for the Conrad
bases $\{G_{E,n}(x)\}$ mentioned above (see also Equation \ref{conradbasis}).

\begin{example}\label{mahlerex}
Let $K=\Qp$, $R_K=\Zp$, $\vert?\vert=\vert?\vert_p$, etc., and let $\{\binom{x}{i}\}$ be the standard Mahler
basis for ${\mathcal C}(\Zp,L)$. By definition, $i!\binom{x}{i}\in \Z[x]$
and is monic of degree $i$; 
thus $\{i!\binom{x}{i}\}$ is an orthonormal basis for
the Tate algebra $L\langle \langle x \rangle \rangle$ of functions regular
on the closed 
unit disc; i.e., locally analytic functions on $\Zp$ of order $0$.
In particular such a function $f$ can then be written
\begin{equation}\label{amice2}
f(x)=\sum_{i\geq 0} c_ii!\binom{x}{i}\,,\end{equation}
where $c_i\to 0$ as $i\to \infty$. By Equation \ref{tatenorm}
we see that if $\{a_i\}$ are the Mahler coefficients of $f$ then
$a_i=i!c_i$. Standard results on the $p$-adic valuation of $i!$ then
 give a simple proof of
Part 2 of Theorem \ref{amice1} in the case $h=0$. The general proof 
in \cite{am1} is given along similar lines.\end{example}

\begin{rem}\label{unifestimates}
Let $f(x)=\sum a_nQ_n(x)$ be locally-analytic of some order $h\geq 0$
as above. Then
part 1 of Theorem \ref{amice1}, and the above example, make it
clear that the estimates on $v(a_n)$ depend {\em only} on $\Vert f\Vert_h$.
Thus if we have a family $\{f_i(x)=\sum a_{n,i} Q_n(x)\}$ of locally analytic
functions of fixed order $h$ {\em and} with
constant (or bounded) norm, then the estimates we obtain on $v(a_{n,i})$ are 
independent of $i$. This observation is essential for our
main result, Theorem \ref{main2}.\end{rem}

Finally, recall that in Equation \ref{intdef} we expressed 
the integral of a continuous function $f$ against a measure $\mu$
as $\sum a_nb_n$ where $f(x)=\sum a_n Q_n(x)$ and $\{b_n\}$ are the
measure coefficients associated to $\{Q_n(x)\}$. The impact of Amice's
Theorem is the following. Let $f$ be locally analytic so that
$a_n\to 0$ very quickly. Then we may integrate all such functions $f$ against
a distribution with coefficients $\{b_n\}$ (defined in the obvious
sense) which may 
{\em not} be bounded, as long as $\sum a_n b_n$ converges. Such distributions
are said to be {\it tempered} and they can readily be described 
(see \cite{am2}). This is
the primary technique used in our main result Theorem \ref{main2}.

\subsection{Appendix: Associated modules}

The results in this appendix elaborate some of the structure associated to
the convolution product of measures. In particular, they explain some earlier
calculations \cite{go8} involving measures in the characteristic $p$ theory.
They are not, however, used in the proof of our main results.

We shall explain here how the 
convolution construction on measures also allows one to make 
${\mathcal C}(R_K,L)$ into a natural ${\mathcal M}(R_K,L)$-module.

\begin{defn}\label{partialconvolution}
Let $f\in {\mathcal C}(R_K,L)$ and let $\mu\in {\mathcal M}(R_K,L)$.
We define $\mu \ast f\in {\mathcal C}(R_K,L)$ by
\begin{equation}\label{partial2}
\mu \ast f(x)=(\mu \ast f)(x):=\int_{R_K} f(x+y)\, d\mu(y)\,.
\end{equation}
\end{defn}

\begin{rem}\label{remark1}
We have used the notation $(\mu,f)\mapsto \mu\ast f$ to distinguish
Definition \ref{partialconvolution}, which associates a continuous
function in $x$ to $(\mu,f)$, from the usual scalar-valued pairing 
$(\mu,f)\mapsto \int_{R_K} f(x)\, d\mu(x)$. Note also that
constructions similar to Definition \ref{partialconvolution} are well known
in classical analysis.
\end{rem}

For instance, Definition \ref{partialconvolution} immediately
gives $\delta_\alpha\ast f(x)=f(x+\alpha)$
for $\alpha \in R_K$ and Dirac measure $\delta_\alpha$. Moreover, for
$\alpha\in R_K$ one sees that
\begin{equation}\label{parametermeasure}
\mu \ast f (\alpha)=\int_{R_K} f(x)\, d\delta_\alpha \ast \mu (x)\,.
\end{equation}

Conversely, of course, the space of measures, ${\mathcal M}(R_K,L)$,
is a natural
${\mathcal C}(R_K,L)$-module where $(f,\mu)\mapsto f(x)\, d\mu(x)$
as usual.

Let $\mu$ be a measure as above and $f\in {\mathcal LA}_h(R_K,L)$.
Using Amice's result, Theorem \ref{amice1}, one sees that
$\mu\ast f$ also belongs to ${\mathcal LA}_h(R_K,L)$.

\subsubsection{The characteristic $p$ case}\label{hyperder}
Definition \ref{partialconvolution} leads to a remarkable differential
formalism in the characteristic $p$ case.
Indeed, let $R_K=\Fr[[T]]$, $K=\Fr((T))$, and $L$ a finite
extension of $K$. The Conrad bases of \cite{co1}
(in particular, the Carlitz polynomials \cite{go8}) lead to the formalism of
{\it differential calculus} in
the above module action. To see this, let
$E=\{e_i\}$ be a Banach basis of the $\Fr$-linear functions
${\mathcal L}_{\Fr}(R_K,L)$. Let $n$ be a non-negative integer
written $r$-adically as $\sum_{t=0}^m b_tr^t$, $0\leq b_t\leq r-1$ all $t$.
Following Carlitz one sets
\begin{equation}\label{conradbasis}
G_{E,n}(x):=\prod_{t=0}^m e_t(x)^{b_t}\,.\end{equation}
Conrad shows that $\{G_{E,n}(x)\}$ is then a Banach basis for 
${\mathcal C}(R_K,L)$.

Standard congruences for binomial coefficients, and the linearity
of $e_i(x)$ all $i$, now immediately imply that
$\{G_{E,n}(x)\}$ satisfies the binomial theorem; that is
$$G_{E,n}(x+y)=\sum_{i=0}^n \binom{n}{i}G_{E,i}(x)G_{E,n-i}(y)\,.$$
As $\{G_{E,n}(x)\}$ is a Banach basis for ${\mathcal C}(R_K,L)$, each
continuous $f$ can be written uniquely as $f(x)=\sum_{n=0}^\infty
b_nG_{E,n}(x)$ where $b_n\to 0$ as $n\to \infty$. Finally we formally set
$$\frac{D_E^i}{i!}=\frac{\frac{d^i}{dE^i}}{i!}$$
to be the measure determined by
\begin{equation}\label{divided measure}
\int_{\Fr[[T]]} f(x)\, d\frac{D_E^i}{i!}(x):=b_i\,.
\end{equation}
Then Equation \ref{partial2} immediately implies
\begin{equation}\label{differentiation}
\frac{D_E^i}{i!}\ast G_{E,n}(x)=\binom{n}{i} G_{E,n-i}(x)\,.
\end{equation}
The analogy with the usual divided-derivatives on power-series is now
obvious.

The reader can easily check that, as operators on 
${\mathcal C}(R_K,L)$, one has
$\displaystyle \frac{D_E^i}{i!}\frac{D_E^j}{j!}=\frac{D_E^{j+i}}{(j+i)!}\,.$
This establishes again 
that the convolution algebra of $L$-valued measures
on $\Fr[[T]]$ 
is isomorphic to the algebra of formal divided derivatives \cite{go8} (which is itself
isomorphic to the algebra of formal divided power series).
Moreover, if $E^\prime$ is another Banach basis for ${\mathcal A}(R_K,L)$,
one may readily express the operator $\frac{D_{E^\prime}^i}{i!}$ in terms of 
$\{\frac{D_E^i}{i!}\}$, and vice versa, by using the co-ordinate free
definition (\ref{partial2}).

\begin{defn}\label{transform1}
Let $z\in R_L$ and let $E$ be as above. Define $\mu_{E,z}$ to be the unique
measure given by 
$$\int_{R_K} G_{E,k}(x)\, d\mu_{E,z}=z^k\,,$$
for non-negative $k$.
If $f(x)=\sum b_n G_{E,n}(x)\in
{\mathcal C}(R_K,L)$, define
\begin{equation}\label{transform2}
\hat{f}_E(z):=\int_{R_K} f(x)\,d\mu_{E,z}(x)=\sum b_n z^n\,.
\end{equation}\end{defn}

The map $f\mapsto \hat{f}_E(z)$ gives a Banach space isomorphism between
${\mathcal C}(R_K,L)$ and the Tate-algebra $L\langle \langle z
\rangle\rangle$ of power-series in
$z$ converging on the unit disc. With this isomorphism,
the operator $\frac{D_E^i}{i!}$ 
transforms into the usual divided-derivative operator
$\frac{D_z^i}{i!}=\frac{\frac{d^i}{dz^i}}{i!}$; i.e.,
\begin{equation}\label{eq1}
\widehat{\frac{D_E^i}{i!}\ast f}(x)=\frac{D_z^i}{i!}\hat{f}(z)\,.\end{equation}
Finally, note that the divided-derivative $\frac{D_z^i}{i!}$ can also
be obtained via digit expansions (as in Equation \ref{conradbasis})
using the Banach basis $\{z^{r^i}\}$ for the space of
$\Fr$-linear elements inside $L\langle \langle z \rangle \rangle$.

\subsubsection{The $p$-adic case} \label{padic}
Let $R_K=\Zp$ and $K=\Qp$, etc. 
Using the Mahler basis $\{\binom{x}{i}\}$,
it is very well-known that the convolution algebra
of $L$-valued measures on $\Zp$ is isomorphic to
${\mathcal R}:=L\otimes_{\Zp} \Zp[[X]]$. Here a measure $\mu$ corresponds to the
power series
$$F_\mu(X)=\sum_{k=0}^\infty\left(\int_{\Zp} 
\binom{x}{k}\,d\mu(x)\right)X^k\,.$$
Clearly, then, we need only compute the action of the measure
associated to $X$ on
${\mathcal C}(\Zp,L)$. Thus let $f(x)=\sum_{k=0}^\infty c_k\binom{x}{k}$
be a continuous function where $\{c_k\}\subset L$ and $c_k\to 0$ as
$k\to \infty$. As $\binom{x+y}{n}=\sum_{j=0}^n\binom{x}{n-j}\binom{y}{j}$,
one immediately computes
\begin{equation}\label{actionX}
X\ast f(x)=\sum_{k=0}^\infty c_k\binom{x}{k-1}=\Delta f(x)\,,
\end{equation}
where $\Delta f(x):=f(x+1)-f(x)$ is the usual difference operator.

Let  $m\in L$ have $\vert m\vert_p<1$ and let $f_m(x):=(1+m)^x$. 
Let $\mu$ be a $\Qp$-valued
measure on $\Zp$ corresponding to a formal power series $F_\mu(X)$. One
checks easily that
\begin{equation}\label{eigenfunc}
\mu\ast f_m(x)=F_\mu(m)f_m(x)\,,\end{equation}
so that the functions $f_m(x)$ are eigenfunctions for the operators
$T_\mu\colon f\mapsto \mu\ast f$. It is simple to see that, up to scalars, 
these are all the common eigenfunctions for
$\{T_\mu\}_{\mu\in {\mathcal M}(\Zp,L)}$ defined over $L$.

\begin{rem}\label{vishik} In Part 4 of the Appendix to \cite{ko1}, there
is an exposition of the $p$-adic spectral theorem of Vishik. Let $\Cp$
be the completion of the algebraic closure of $\Qp$ with the canonical
extension of $\vert?\vert_p$. Vishik's theorem applies to ``analytic
operators;'' i.e., operators $A$ with compact spectrum $\sigma_A$
over $\Cp$ and with analytic resolvent $R_A(z):=(z-A)^{-1}$ on the complement
of $\sigma_A$. (Here analyticity means essentially 
that the matrix elements of the operator are Krasner analytic; for more
see \cite{ko1}.) Let $\mu$ be the measure on $\Zp$ corresponding to
$X$ as above. One checks that the spectrum of the operator
$f\mapsto \mu\ast f=\Delta f$ over $\Cp$ is given
by the eigenvalues
$\{x\in \Cp\mid \vert x\vert_p<1\}$ which is obviously not a
compact set. It is however, bounded and therefore
it is reasonable to ask whether a form of
Vishik's results might also hold for $\Delta$.
\end{rem}

One can also express the action of $L\otimes\Zp[[X]]$ on ${\mathcal C}(\Zp,L)$ via
the following construction. Let $B$ be the Banach space of bounded
sequences $\{b_i\}_{i\in \Z}\subset L$ equipped with the sup norm.
We write these sequences formally as
$f(X,X^{-1})=\sum_{i=-\infty}^\infty b_i X^i$. 

Let $H$ be the subspace of all $f(X,X^{-1})=\sum b_iX^i\in B$ with
$b_i\to 0$ as $i\to -\infty$ (note the minus sign!). In other words,
$H$ consists of all $f(X,X^{-1})$ whose polar part converges for
$X^{-1}\in R_L$. A little thought establishes that $H$ is actually
a closed $\mathcal R$-submodule of $B$ where the action
of ${\mathcal R}:=L\otimes\Zp[[X]]$ is via multiplication of power series in the usual sense.
Thus $H/X\mathcal R$ is
Banach space isomorphic to the Tate algebra $L\langle \langle X^{-1}\rangle
\rangle$ and equips $L\langle \langle X^{-1}\rangle \rangle$ with a 
natural $\mathcal R$-module structure. Intuitively we simply
multiply the two series in the ``usual'' fashion and throw out the terms $X^i$ with $i$ positive.

 Let $X^{-1}\in R_L$ and let 
$f(x)=\sum_{k=0}^\infty c_k\binom{x}{k}\in {\mathcal C}(R_K,L)$.
Let $\mu_{_X}$ be the measure on $\Zp$ given by
$$\int_{\Zp} \binom{x}{k}\, d\mu_{_X}(x)=X^{-k}\,,$$
for non-negative $k$.
Set $\hat{f}(X)=\sum_{k=0}^\infty c_kX^{-k}=\int_{\Zp} f(x)\, d\mu_{_X}(x)$, 
as in Equation \ref{transform2},
which again gives a Banach space isomorphism between the continuous
functions and the Tate algebra in $X^{-1}$. With this isomorphism, 
the action of
the measure associated to $F(X)=\sum b_iX^i$ on $f$ transforms into the 
action of $F(X)\in \mathcal R$ on $\hat{f}(X)$ presented above. Our thanks
to W.\ Sinnott for pointing this construction out to us.

\subsubsection{A curious connection between continuity on
$\Zp$ and finite characteristic measures}\label{curious}
Let $L$ be a finite extension of $\Qp$ and let $f\colon \Zp\to L$
be a continuous function written $\sum b_n \binom{x}{n}$. Set
$a_n:=f(n)$ for $n=0,1...$ and form two divided power
series $A(x):=\sum a_n\frac{x^n}{n!}$ and $B(x):=\sum b_n\frac{x^n}{n!}$.
It is well-known that $e^{-x}A(x)=B(x)$. Conversely, if one defines
$B(x)$ by this equation for any sequence $\{a_n\}$, then we obtain
a continuous function $f$ with $f(n)=a_n$ if and only if
$b_n$ tends to $0$ as $n\to \infty$.

Such a formalism also works when $L$ has
finite characteristic (see \S 8.4 of \cite{go4}). 
The results, and notation,
 of Subsection \ref{hyperder} then give the following curious
result. Let $\mu$ be a measure on $R_K$ and let $\{b_n\}$ be its measure
coefficients with respect to $\{G_{E,n}(x)\}$. Let $\nu$ be the unique
measure with coefficients $\{(-1)^n\}$ with respect to $\{G_{E,n}(x)\}$.
Then {\it there is a continuous $L$-valued function $f\colon \Zp \to L$ with
$f(n)=b_n$, all $n$,
if and only if the measure coefficients of $\nu\ast \mu$ tend to
$0$.}

The condition that the measure coefficients tend to $0$ is stronger than what
is needed to integrate continuous functions. Perhaps there is a larger
class of functions that may be integrated by such a measure.

\section{$L$-functions of $\tau$-sheaves}\label{ltau}
\subsection{$\tau$-sheaves}\label{deftausheaf}
The concept of a $\tau$-sheaf (\cite{tw1}, \cite{tw2}, \cite{bp1}, \cite{boc1},
\cite{ga1}, \cite{ga3}, \cite{go8})
arose out of the concept of a Drinfeld 
module which is where we begin. Let $C$ be a smooth, projective, geometrically
connected curve over the finite field $\Fr$ where $r=p^{m_0}$ with
$p$ prime. As usual one chooses a fixed closed point $\infty$ of
degree (over $\Fr$) $d_\infty$. The space $C^\prime:=C-\infty$ is affine 
and one denotes by $A$ the ring of functions which are regular on $C^\prime$.
As is well known, the ring $A$ is a Dedekind domain with 
unit group $\Fr^\ast$ and
finite class group. Let $k$ be the quotient field of $A$.

A field $L$ with an $\Fr$-algebra map $\imath\colon A\to L$ is said to be
an ``$A$-field;'' the kernel $\mathfrak p$ of $\imath$ is a prime ideal of
$A$ which is called the ``characteristic of $L$.'' 
Let $\bar{L}$ be a fixed algebraic closure of $L$ (which
is obviously also an $A$-field with the same characteristic) 
and let $\tau\colon \bar{L}\to \bar{L}$ be
the $r$-th power mapping, $\tau (l)=l^r$. The elements $l\in L$ and $\tau$
generate, by composition, an algebra of endomorphisms $L\{\tau\}$ 
of $\bar{L}$ with
$\tau l=l^r\tau$, etc. There is unique homomorphism $D\colon L\{\tau\}
\to L$ given by $D(\sum_{j=0}^t b_j \tau^j)=b_0$. A {\it Drinfeld $A$-module}
$\psi$ over $L$ is an injection of $\Fr$-algebras $\psi \colon A\to
L\{\tau\}$, $a\mapsto \psi_a(\tau)$, such that $D\circ \psi=\imath$ 
but $\psi_a\neq \imath(a) \tau^0$ for some $a\in A$. 
We denote by $\psi[a]$ the finite subgroup of elements $z \in \bar{L}$ with
$\psi_a(z)=0$. As $A$ is obviously commutative, $\psi[a]$ inherits
an $A$-module structure. One can show the
existence of an integer $d>0$ such that $\psi[a]$ is 
$A$-module isomorphic to $A/(a)^d$ for all $a\in A-\mathfrak p$. 
One calls $d$ {\it the rank of $\psi$.}

The next key step was taken by Drinfeld and then 
Anderson. Let $M:=L\{\tau\}$ which we
now view as an $L\otimes_\Fr A$-module in the following fashion. Let
$f(\tau)=\sum_{i=0}^j c_i\tau^i$ be an arbitrary element of $M$, $a\in A$,
and $l\in L$. One sets
\begin{equation}\label{definingstruc}
l\otimes a \cdot f(\tau):=lf(\psi_a(\tau))\,.\end{equation}
One checks that $M$ then becomes a projective $L\otimes A$-module of
rank $d$ and thus gives rise to a locally-free sheaf on $\text{
Spec}(L\otimes A)$. See \cite{an1}, for example, 
for more on how the properties of $\psi$ may
be reinterpreted in terms of $M$. 

The module $M$ also possesses an obvious action by $\tau$, with
$(\tau, f(\tau))\mapsto \tau f(\tau)$ (multiplication in $L\{\tau\}$). Note
that for $l\in L$ one has
\begin{equation}\label{tau1}
\tau l f(\tau)=l^r f(\tau)\,.\end{equation}

The essential features of $M$ is that it is a coherent sheaf with
``$A$-coefficients'' equipped with a $\tau$-action as above. 
This leads directly
to the basic notion of a $\tau$-sheaf. Thus let $X$ be a scheme over
$\Fr$. Let $\sigma=\sigma_X$ be the absolute
Frobenius morphism with respect to $\Fr$; that is for any affine 
$\text{Spec}(B)\subseteq X$ and $b\in B$ one has $\sigma^\ast b=b^r$.
\begin{defn}\label{tausheaf} 
A {\it $\tau$-sheaf} on $X$ is a pair $\underline{\mathcal F}:=
({\mathcal F},\tau)$ consisting of a locally-free sheaf $\mathcal F$ on
$X\times_\Fr C^\prime$ and an $({\mathcal O}_X\otimes A)$-linear morphism
\begin{equation}\label{tauaction2}
\tau=\tau_{_{\mathcal F}}\colon (\sigma \times \text{id})^\ast
{\mathcal F}\to \mathcal F\,.\end{equation} 
A {\it morphism of $\tau$-sheaves} is a morphism of the underlying
coherent sheaves which commutes with the $\tau$-actions.
\end{defn}
\noindent
The reader will readily see that the Frobenius-linear property
expressed in Equation \ref{tau1} is equivalent to the formulation
given in Equation \ref{tauaction2}.

\begin{rem}\label{strict}
In the papers \cite{boc1}, \cite{boc2},\cite{bp1} and \cite{go7} a more
general notion of $\tau$-sheaf is given where the underlying module
need only be a coherent module. 
\end{rem}

\subsection{Domain spaces}\label{domain}
We will present here the basic ideas on exponentiation of ideals as in
Section 8.2 of \cite{go4} and \cite{boc1}. 

\subsubsection{The theory at $\infty$}\label{domaininfty}
We begin with the place $\infty$. Let $K:=k_\infty$ now denote 
the completion of $k$ at $\infty$ 
and let $\pi\in K$ now denote a fixed uniformizer of $K$; this is notation
that henceforth will be used throughout the rest of this paper. 
Let $\F_\infty=\F_K\simeq \F_{r^{d_\infty}}$ be the associated finite field. Thus every element
$\alpha\in K^\ast$ can be written uniquely as
\begin{equation}\label{alpha}
\alpha=\zeta_\alpha \pi^{n_\alpha}\langle \alpha \rangle\,,\end{equation}
where $\zeta_\alpha\in \F_\infty^\ast$, $n_\alpha=v_\infty(\alpha) \in \Z$, and
$\langle \alpha \rangle$ belongs to the group $U_1(K)$ of $1$-units of $K$. 
Note that both $\zeta_\alpha$ and $\langle \alpha \rangle$ depend on the
choice of $\pi$. For $\alpha \in k^\ast$ we set
\begin{equation}\label{deg}
\deg_k (\alpha)=-d_\infty n_\alpha\,;\end{equation}
as usual we set $\deg_k(0)=-\infty$.

For any non-zero fractional ideal
 $I$ of $A$, we also let $\deg_k (I)$ be the degree over $\Fr$ of
the divisor associated to $I$ on the curve $C^\prime$; 
thus $\deg _k(I)=\deg_{\Fr} A/I$
when $I\subset A$. Moreover
Equation \ref{deg} implies that for $\alpha\in k^\ast$ one has
\begin{equation}\label{deg2}
\deg_k (\alpha)=\deg_k (\alpha A)\,.\end{equation}
Thus $\deg_k (\alpha)$ is the degree of the finite part of the divisor
of $\alpha$ on the complete curve $C$.

We let $\C_\infty$ be the completion of a fixed algebraic closure $\bar{K}$
equipped with the canonical extension of the normalized
absolute value $\vert?\vert_\infty$.

The element $\alpha$ is said to be {\it positive} (or {\it monic}) if
and only if $\zeta_\alpha=1$ (so that the notion of positivity most definitely
depends on the choice of $\pi$). Clearly the product of two
positive elements is also positive. 

Let $\alpha$ be positive.

\begin{defn}\label{sinfty}
1.\ We set $S_\infty:=\C_\infty^\ast \times \Zp$\,.\\
2.\ Let $s=(x,y)\in S_\infty$. We then set
$$\alpha^s:=s^{\deg_k (\alpha)}\langle \alpha \rangle^y\,.$$
\end{defn}

\noindent
As $\langle \alpha \rangle=1 + \lambda_\alpha$ with $v_\infty (\lambda_\alpha)
>0$, one has the convergent expression
\begin{equation}\label{lambdaalpha}
\langle \alpha \rangle^y=\sum_{j=0}^\infty \binom{y}{j}\lambda_\alpha^j\,.
\end{equation}

The space $S_\infty$ will be the domain for the
$L$-series of $\tau$-sheaves at $\infty$.

The group action on $S_\infty$ will be written additively. Suppose
that $j\in \Z$ and $\alpha^j$ is defined in the usual sense of the
canonical $\Z$-action on the multiplicative group. Let
$\pi_\ast\in \C_\infty^\ast$ be a fixed $d_\infty$-th root of $\pi$. 
Set $s_j:=(\pi_\ast^{-j},j)\in S_\infty$. One checks easily that 
Definition \ref{sinfty} gives $\displaystyle \alpha^{s_j}=\alpha^j$. 
When there is no chance of confusion, we denote $s_j$ simply by ``$j$.''

Let $\mathcal I$ be the group of fractional ideals of the Dedekind
domain $A$ and let ${\mathcal P}\subseteq \mathcal I$ be the subgroup of
principal ideals. Let ${\mathcal P}^+\subseteq \mathcal P$ be the
subgroup of principal ideals which have positive generators. One knows
that ${\mathcal I}/{\mathcal P}^+$ is a finite abelian group. The 
association 
\begin{equation}\label{langle1}
{\mathfrak h}\in {\mathcal P}^+\mapsto \langle {\mathfrak h}\rangle:= 
\langle \lambda \rangle\,,\end{equation}
 where $\lambda$ is the unique positive generator of $\mathfrak h$,
is obviously a homomorphism from ${\mathcal P}^+$ to $U_1(K)$. 

For the moment, let $u=1+m \in U_1(\C_\infty)$, $\vert m\vert <1$,
 be an arbitrary $1$-unit
in $\C_\infty$ and let $y=\sum_{j=j_0}^\infty c_j p^j$ be an arbitrary
element of $\Qp$. One sets $u^y:=\prod_j (1+m^{p^j})^{c_j}$ which obviously
converges in $\C_\infty$. Thus $U_1(\C_\infty)$ is naturally a 
$\Qp$-vector space; in particular, $U_1(\C_\infty)$ is thereby
a divisible, and thus injective, group. We therefore have the next result.

\begin{lemma}\label{langle2}
The homomorphism $\langle ?\rangle\colon {\mathcal P}^+\to U_1(\C_\infty)$
extends uniquely to a homomorphism $\langle ?\rangle\colon
{\mathcal I}\to U_1(\C_\infty)$.\end{lemma}

\noindent
The uniqueness in Lemma \ref{langle2} follows from the fact that
${\mathcal P}^+$ has finite index in $\mathcal I$.

\begin{defn}\label{sinfty2}
Let $I\in {\mathcal I}$ and $s=(x,y)\in S_\infty$. We then set
\begin{equation}\label{sinfty3}
I^s:=x^{\deg_{k}I}\langle I\rangle^y\,.\end{equation}\end{defn}

The reader will easily see that the mapping
$${\mathcal I}\times S_\infty\mapsto I^s$$ is bilinear.

\begin{defn} \label{valuefield}
Let ${\mathbb V}\subset \C_\infty:=k(\{I^{s_1}\mid I\in {\mathcal I}\})$. We
call $\mathbb V$ the {\it value field} associated to $\pi$ and $\pi_\ast$.
The place on $\mathbb V$ given by its inclusion in $\C_\infty$ will be 
also be denoted $\infty$ and is called the {\it canonical infinite place of}
$\mathbb V$.
\end{defn}

\begin{prop}\label{value1}
The field $\mathbb V$ is a finite extension of $k$.
\end{prop}

\begin{proof}
If $I=(i)$ where $i$ is positive, then $I^{s_1}=i\in k$. As
${\mathcal I}/{\mathcal P}^+$ is finite, the result follows. \end{proof}

Let ${\mathcal O}_{\mathbb V}\subset \mathbb V$ be the ring of
$A$-integers. The places of ${\mathbb V}$ which lie outside of
$\text{Spec}({\mathcal O}_{\mathbb V})$ (and so lie over the place
$\infty$ of $k$) are the ``infinite primes of ${\mathbb V}$;''
thus places lying above $\text{Spec}(A)$ are the ``finite primes.''

Let $\alpha $ be an element of $\mathbb V$. We
let $\deg_{\mathbb V}(\alpha)$ be the degree over
$\Fr$ of the finite part of the
divisor of $\alpha$; as the degree of a principal divisor is $0$, 
this is the opposite of the degree of the infinite part of
the divisor of $\alpha$. In particular, for $\alpha \in k$, one has
\begin{equation}\label{degree1}
\deg_{\mathbb V}(\alpha)=[{\mathbb V}\colon k]\deg_k(\alpha)\,.
\end{equation}
Similarly, if $J$ is a fractional ${\mathcal O}_{\mathbb V}$-ideal, then
we let $\deg_{\mathbb V}(J)$ be the degree over
$\Fr$ of the associated finite divisor,
etc. If $I$ is an $A$-fractional ideal then we again have
\begin{equation}\label{degree2}
\deg_{\mathbb V}(I{\mathcal O}_\mathbb V)=[{\mathbb V}\colon k]\deg_k(I)\,.
\end{equation}

The next proposition and corollary, which are elementary, are explicitly
recorded as they will be used in the proof of our main result.

\begin{prop}\label{nonzero}
Let $0\neq\alpha\in {\mathcal O}_{\mathbb V}$. Then
$\deg_{\mathbb V}(\alpha)\geq 0$.
\end{prop}

\begin{cor}\label{nonzero2}
The only element of ${\mathcal O}_\mathbb V$ of negative degree
is $0$.
\end{cor}

\begin{prop}\label{value2}
Let $I$ be a non-trivial ideal of $A$.
Then, 
\begin{equation} \label{value3}
I{\mathcal O}_{\mathbb V}=I^{s_1}{\mathcal O}_\mathbb V\,.\end{equation}
Moreover, $I^{s_1}$ has a pole at every infinite place of $\mathbb V$.
Finally, let $\mathfrak P$ be a prime of ${\mathcal O}_{\mathbb V}$
with additive valuation $v_{\mathfrak P}(?)$.
Then 
\begin{equation}\label{value4}
\deg_k({\mathfrak P})
v_{\mathfrak P}(I^{s_1})\leq [{\mathbb V}\colon k]\deg_{k}(I)\,.\end{equation}
\end{prop}

\begin{proof} Let $t$ be the order of ${\mathcal I}/{\mathcal P}^+$;
thus $I^t=(\alpha)$ where $\alpha\in A$ is positive. By definition
one has $\lambda^t=\alpha$ where $\lambda=I^{s_1}$. Thus the order
of $I^{s_1}$ is the same as the order of $I$ at each prime of
${\mathcal O}_{\mathbb V}$, and Equation \ref{value3} follows. Clearly
$\alpha$ has a pole at every infinite prime of $\mathbb V$ and so
therefore must $I^{s_1}$. Finally,
Equation \ref{value3} implies that $\deg_{\mathbb V}(I{\mathcal O}_\mathbb V)
=\deg_{\mathbb V}(I^{s_1})$ as a principal divisor has degree $0$. Thus Equation \ref{value4} follows 
from Equation \ref{degree2}.\end{proof}

In particular, while the elements $\{I^{s_1}\}$ are not necessarily
in $A$, they do behave very much like elements of $A$. For instance,
they have the same absolute value at the (non-normalized) extension
of $\vert ?\vert_\infty$ at {\em any} infinite place of $\mathbb V$. This will
be of importance in the proof of our main result.

\subsubsection{The theory for finite primes}\label{domainfinite}
Let $v$ be the place associated to a prime $\mathfrak p$ of $A$ and
set $d_v:=\deg_{\Fr}(v)$. 
Let $k_v$ be the associated completion of $k$ with normalized
absolute value $\vert ?\vert_v$. Let $\bar{k}_v$ be a fixed algebraic
closure of $k_v$ and let ${\mathbb C}_v$ be its completion with the
canonical extension of $\vert ?\vert_v$. Finally,
let $\sigma\colon {\mathbb V}\to \bar{k}_v$ be an embedding over $k$ and
set 
\begin{equation}\label{ksigmav}
k_\sigma=k_{\sigma,v}:=k_v(\sigma ({\mathbb V}))\,.\end{equation}
By Proposition \ref{value1} one sees that $k_{\sigma,v}$ is finite
over $k_v$ and one lets $A_\sigma=A_{\sigma, v}\subset k_{\sigma,v}$ 
denote the ring of $A_v$-integers, with maximal ideal
$M_{\sigma,v}$, and residue degree
$f_\sigma=f_{\sigma, v}$. Any element $\beta\in A_{\sigma,v}^\ast$
can then be written 
\begin{equation}\label{decompbeta}
\beta=\omega_{\sigma,v}(\beta)\langle \beta\rangle_{\sigma,v}\,,
\end{equation}
where $\omega_{\sigma,v}(\beta)$ belongs to the group $\mu_{r^{d_vf_\sigma}-1}$
of roots of unity inside $A_{\sigma,v}$, and where 
$\langle \beta \rangle_{\sigma,v}$ is a $1$-unit.

\begin{defn}\label{ssubv}
We set
\begin{equation}\label{ssubv2}
S_\sigma=S_{\sigma,v}:=\varprojlim_j \Z/\left(p^j(r^{d_vf_\sigma}-1)\right)\simeq \Zp\times
\Z/(r^{d_vf_\sigma}-1)\,.\end{equation}
\end{defn}

Let $y_\sigma=(y_{\sigma,0},y_{\sigma,1})\in S_\sigma$ and let $\beta$ be as in
Equation \ref{decompbeta}. Then we set
\begin{equation}\label{vadicexp}
\beta^{y_\sigma}:=\omega_{\sigma,v}(\beta)^{y_{\sigma,1}}\langle \beta
\rangle_{\sigma,v}^{y_{\sigma,0}}\,.\end{equation}
Let ${\mathcal I}(v)$ be the group of $A$-fractional ideals generated
by the primes $\neq v$ and let $I\in {\mathcal I}(v)$. One knows
that $\sigma(I^{s_1})\in A_{\sigma,v}^\ast$ by Equation \ref{value3}.
Let $s_\sigma=(x_v,y_\sigma)\in \C_v^\ast\times S_\sigma$. Finally we define
\begin{equation}\label{vadicexp2}
I^{s_\sigma}:=x_v^{\deg_{\Fr}(I)}\sigma(I^{s_1})^{y_\sigma}\,.\end{equation}
The space $\C_v^\ast \times S_\sigma$ is the domain for the $v$-adic
theory of $L$-series associated to $\tau$-sheaves.

Note that if $I=(i)$ with $i$ a positive $v$-adic unit, then
$$I^{s_\sigma}=x_v^{\deg_{k}(i)}i^{y_\sigma}\,.$$
In particular, one has $I^{(1,j,j)}=i^j$ for all integers $j$.

\subsubsection{Entire functions}\label{entire}
We begin here with the theory for $S_\infty$. The theory
for $\C_v^\ast\times S_\sigma$ is entirely similar and will be left
to the reader. The basic reference is Section 8.5 of \cite{go4} (but
see also \cite{boc1} in this regard).

Our $L$-series at $\infty$
will be ``Dirichlet series'' of the form $\sum c_II^{-s}$ 
where $s=(x,y)\in
S_\infty=\C_\infty^\ast \times \Zp$. By Definition \ref{sinfty2}, 
$I^{-s}=x^{-\deg_k I}\langle I \rangle^{-y}$; therefore Dirichlet
series immediately give rise to a $1$-parameter family of formal
power series in $x^{-1}$. This leads to the following definition.

\begin{defn}\label{entirefunction}
An {\it entire function on $S_\infty$} is a $\C_\infty$-valued uniformly
convergent family of entire power series in $x^{-1}$ parameterized
by $\Zp$.
\end{defn}

\noindent
As explained in \cite{go4} this means that for each $y\in\Zp$ we are
given a power series $g_y(1/x)$ where $g_y(u)$ converges for all
$u$. Moreover, for each bounded set $B\subset \C_\infty$ and
$\epsilon>0$ there exist
a $\delta:=\delta_B>0$ such that if $y_0$ and $y_1$ are in $\Zp$
with $\vert y_0-y_1\vert_p<\delta$ then $\vert g_{y_0}(u)-g_{y_1}(u)\vert_\infty<\epsilon$ for all $u\in B$. This forces the zeroes to ``flow continuously.''

Let $0\neq \rho\in \C_\infty$. Then, exactly as in Definition \ref{analB},
one can define a norm on the space of functions analytic on 
the closed disc $\{u\in \C_\infty\mid \vert u\vert_\infty\leq 
t=\vert \rho\vert\}$. In
\cite{go4}, Definition \ref{entirefunction} is reinterpreted in terms of 
the family of norms defined by all such $\rho$.

In \cite{boc1}, one chooses
a family $\{\rho_j\}_{j=0}^\infty \subset \C_\infty^\ast$ of increasing
and unbounded norm to create a Fr\'echet space out of the entire
power series in $u$. An entire function on $S_\infty$, as in Definition \ref{entirefunction},
is then just a continuous function from $\Zp$ into this Fr\'echet
space.

As will be seen, the Dirichlet series of $\tau$-sheaves are entire
in the sense of Definition \ref{entirefunction}. Suppose for the
moment that $A=\Fr[T]$ and let 
$\displaystyle L(s)=\sum_{a\in A~\text{positive}}
c_aa^{-s}$, where $\{c_a\}\subset A$, be entire
as above. Let $j$ be a non-negative integer and set
\begin{equation}\label{zfunc}
z_L(x,-j):=L(\pi_\ast^{j}x,-j)\,.\end{equation}
It is easy to see that $z_L(x,-j)$ is now a power series in
$x^{-1}$ with $\Fr[T]$-coefficients. Moreover, as $L(s)$ is
entire in $s$, we conclude that $z_L(x,-j)$ is entire in $x^{-1}$, and,
in particular, its coefficients must tend to $0$ in $\Fr((1/T))$. 
As the degrees of non-zero elements in $A$ are obviously non-negative, we see
that the only way this can happen is that almost all such coefficients
vanish; i.e., $z_L(x,-j)$ must belong to $A[x^{-1}]$. 

In fact, the rationality result just established for $\Fr[T]$ also
holds for all $L$-series of $\tau$-sheaves, for {\em any}
$A$, and leads to the next definition. Let $A$ now be completely general.

\begin{defn}\label{essalgfunction}
Let $L(s)=\sum_I c_I I^{-s}$ be an entire function where 
$\{c_I\}\subset \C_\infty$ lies in a finite extension $H$ of $k$. Put
$H_1:=H\cdot \mathbb V$ and define $z_L(x,-j)$ for $j\geq 0$ exactly as in
Equation \ref{zfunc}. We then say that $L(s)$ is {\it essentially algebraic}
if and only if $z_L(x,-j)\in H_1[x^{-1}]$ for all $j$.
\end{defn}

Let $v$ be a finite prime of $A$. We leave to the reader 
the easy translation of
``entire'' and ``essentially algebraic'' to {\it $v$-adic Dirichlet
series} of the form $\sum_{I\in {\mathcal I}(v)}c_II^{-s_\sigma}$, 
$s_\sigma\in \C_v^\ast\times
S_\sigma$ as in Equation \ref{vadicexp2}.

The idea behind all of this is that one starts with a Dirichlet series 
$L(s)=\sum_I c_I I^{-s}$,
$s\in S_\infty$, with coefficients
in some finite extension $H$ of $k$ and then using the injections $\sigma$
of $H$ into $\bar{k}_v$, one defines the various ``interpolations'' 
$L(s_\sigma)$
of this $L$-series at the finite primes via the {\em same}
sum {\em except} that we have removed the factors lying over $v$;
i.e., $L(s_\sigma):=\sum_{I\in {\mathcal I}(v)} \sigma(c_I)I^{-s_\sigma}$.
Therefore these interpolations are defined
a priori simply as certain $v$-adic Dirichlet series.
(See also Remark \ref{infiniteplaces} just below.)
In general there is no reason to expect any relationship between
the interpolations of $L(s)$ at different places.
However, the $L$-series of a $\tau$-sheaf will turn out to be
an essentially algebraic entire function.  The {\it special polynomials}
$\{z_L(x,-j)\}_{j=0}^\infty$ then have two basic attributes. First of
all, by the work of B\"ockle and Pink\cite{bp1}, one can express them 
{\it cohomologically} (this is how one
deduces in general that the power series $z_L(x,-j)$ is indeed
a polynomial). Secondly, they allow one
to relate the $\infty$-adic and $v$-adic theories by simply removing
the finite number of Euler-factors lying above $v$ in the 
special polynomials and then substituting
$x_v$ for $x$, etc. In other words, one 
obtains the functions $L(s)$ and $L(s_\sigma)$ associated to an
essentially algebraic Dirichlet
series by simply interpolating the special polynomials
$\{z_L(x,-j)\}_{j=0}^\infty$; as the non-positive integers are
{\em dense} in $\Zp$ and $S_v$.
This is obviously not possible more generally.

\begin{rem}\label{infiniteplaces}
In the case where $A$ is a principal ideal domain (e.g., $A=\Fr[T]$),
the collection of interpolations given above ranges over {\em all}
the places of the quotient field $k$. However, when $A$ is general, there
is no reason to believe that $\mathbb V$ has only one infinite place
(corresponding to its inclusion in $\C_\infty$). Thus let
$\sigma\colon {\mathbb V}\to \C_\infty$ be an injection corresponding
to a different infinite place. Let $I$ be an ideal of $A$ and set
\begin{equation}\label{otherinfiniteplace}
\langle I\rangle_\sigma:=\pi_\ast^{\deg_k I}\sigma(I^{s_1})\,.
\end{equation}
Using the fact that ${\mathcal I}/{\mathcal P}^+$ is finite, one concludes
as before that $\langle I\rangle_\sigma$ is a unit (i.e., has norm $1$)
in $\C_\infty$. However, one cannot conclude that $\langle I\rangle_\sigma$
is a $1$-unit; indeed there may be a non-trivial group
of roots of unity
one needs to handle. Nevertheless, these may be handled in exactly the same
fashion as in the $v$-adic theory and one readily defines 
$L$-series associated to the infinite embedding $\sigma$.
We will denote these by $L(\underline{\mathcal F}, s_{\sigma,\infty})$
where $s_{\sigma,\infty}=(x,y_{\sigma,\infty})\in 
\C_\infty^\ast\times\varprojlim_j \Z/(r^{d_\sigma}-1)p^j$ where $d_\sigma$ is the degree
over $\Fr$ corresponding to the place of $\mathbb V$
given by $\sigma$. We will not stress
these functions here as they have no classical counterparts
and they may unnecessarily confuse the reader. However, when they
arise from the $L$-series of a $\tau$-sheaf (as presented 
in the next subsection),
the reader may easily check that our techniques
show that these functions also possess the same features as all other
interpolations.\end{rem}

\subsection{Euler factors}\label{eulerfact}
We now follow \cite{bp1}, \cite{boc1}
and define the $L$-series of a $\tau$-sheaf via
an Euler product. Let $X$ be a scheme of finite type over 
$\text{Spec}(A)$ and let
$\underline{\mathcal F}$ be a $\tau$-sheaf as in Definition \ref{tausheaf}.
Let $X^0$ be the set of closed points of $X$ and for each $\alpha\in X^0$, let
$\mathfrak{p}_\alpha$ be its image in $\text{Spec}(A)$.  

\begin{defn}\label{boeckle1}
We define
\begin{equation}\label{boeckle2}
L(\alpha,\underline{\mathcal F},u)^{-1}:=\det_k (I_d-u\tau\mid {\mathcal F}_\alpha
\otimes_A k)\in k[u]\,,\end{equation}
where ${\mathcal F}_\alpha$ is the fiber of the sheaf $\mathcal F$, $I_d$ is the
identity morphism, and the determinant is taken over $k$.\end{defn}

In \cite{bp1} it is shown that $L(\alpha,\underline{\mathcal F},u)^{-1}\in
A[u^{d_\alpha}]$ where $d_\alpha$ is the degree of $\alpha$ over $\Fr$. Let
$d_{{\mathfrak p}_\alpha}$ be the degree of $\mathfrak p_\alpha$ over $\Fr$.
Notice that $A[u^{d_\alpha}]\subseteq A[u^{d_{{\mathfrak p}_\alpha}}]$

\begin{defn}\label{global}
Let $s\in S_\infty$. We set 
\begin{equation}\label{global2}
L(\underline{\mathcal F},s):=\prod_{\alpha\in X^0}
L(\alpha,\underline{\mathcal F},u)
\vert_{u^{d_{{\mathfrak p}_\alpha}}={\mathfrak p}_\alpha^{-s}}\,.\end{equation}
\end{defn}

B\"ockle shows that the Euler product (\ref{global2}) converges on
a ``half-plane'' of $S_\infty$. That is, there exists a non-zero number 
$t \in \R$
such that (\ref{global2}) converges for all $s\in (x,y)\in S_\infty$ such
that $\vert x\vert_\infty\geq t$. 

\begin{rem}\label{addversion}
For future use, we rewrite the condition for the half-plane in terms of
additive valuations. Let $K$ be the completion of $k$ at
$\infty$, as usual, and put $K_1:=K[\langle I\rangle]$ where
$I$ runs over the ideals of $A$; as ${\mathcal I}/{\mathcal P}^+$ is finite,
$K_1$ is finite over $K$.
(The reader will be tempted to conclude that $K_1$ is the
compositum of $K$ and $\mathbb V$. This is only obviously so when
$d_\infty=1$; the general relationship is not yet clear.) By abuse of notation let
$A_\infty$ denote the maximal compact subring of $K$ and $A_{1,\infty}$
that of $K_1$. Let $e$ be the ramification degree and $v_\infty$ (resp.\
$v_{1,\infty}$) the canonical additive valuation on $K$ (resp.\ $K_1$)
which assigns $1$ to a uniformizing parameter. Thus,
upon extending these canonically to the algebraic closure,
one has  $v_{1,\infty}=ev_{\infty}$.
Recall that we set $q_K=r^{\deg \infty}$. One then sees readily
that 
$$\vert x\vert_\infty\geq t \Leftrightarrow v_{1,\infty}(x)\leq -e\log_{q_K}(t)\,.$$
\end{rem}

Let $v$ now be a finite prime of $A$. The $v$-adic
version of $L(\underline{\mathcal F},s)$, which will be denoted
$L(\underline{\mathcal F},s_\sigma)$, $s_\sigma\in \C_v^\ast\times S_\sigma$,
is now obvious using the embedding $\sigma$. Indeed,
one simply uses in Definition \ref{global},
$\alpha\in X^0(v)$ where
$X^0(v)$ is the set of closed points not lying over $v$. Also obvious
in this case is the existence of a $v$-adic half-plane of convergence as the
Euler factors all have coefficients in ${\mathcal O}_\mathbb V$.

\begin{rem}\label{absconv}
In both the $\infty$-adic and $v$-adic theories, we abuse
language and say that the 
Dirichlet series {\it converges absolutely} in the half-plane of convergence
of its Euler product. Let
$s\in S_\infty$, where $\vert x\vert_\infty\geq t$
($t$ as above), and expand $L(\underline{\mathcal F},s)$ as
$\sum_I b_I I^{-s}$. Then as $\deg I\to \infty$ one has $b_II^{-s}\to 0$
which is what we shall mean by ``absolute convergence.''
An exactly similar statement holds $v$-adically.\end{rem}

The next result is a restatement of one of the main theorems of \cite{boc1}
and we refer the reader there for the proof.

\begin{theorem}\label{cohenmc1}
Let $X$ be a reduced, affine, equi-dimensional Cohen-Macaulay variety over
$\text{Spec}(A)$ of dimension $e_X$. Let $\underline{\mathcal F}$ be a 
$\tau$-sheaf (which we recall is locally-free by definition in this
paper). Then both
$L(\underline{\mathcal F},s)^{(-1)^{e_X-1}}$ and 
$L(\underline{\mathcal F},s_\sigma)^{(-1)^{e_X-1}}$ are essentially algebraic
entire functions. Moreover, for each non-negative integer $j$, the 
degree in $x^{-1}$ (respectively $x_v^{-1}$) of the associated
special polynomial at $-j$ is $O(\log (j))$.
\end{theorem}

\begin{rems}\label{remarks1}
1.\ In \cite{boc1} a slightly more restricted choice of local
uniformizer $\pi$ is chosen. This allows for certain global rationality
statements that we have not given here. In any case, the arguments in
\cite{boc1}, and in particular, the growth estimate of the special
polynomials, actually apply in complete generality as we have set things up.\\
2. The definition of the $L$-factors of $\tau$-sheaves given in 
Definition \ref{boeckle1} seems very different than the usual definition
of, say, the $L$-factors of elliptic curves where one uses 
Tate modules and Frobenius morphisms etc. For $\tau$-sheaves one
can also use this approach and indeed one obtains the same local
factors \cite{boc1}. For instance, the $L$-series of a Drinfeld module $\psi$
has traditionally
been defined this way. At the good primes, it is relatively easy to see
that both approaches agree. At the bad primes, one needs to use the 
{\it Gardeyn maximal model} (e.g., \cite{ga1})
of the $\tau$-sheaf associated $\psi$ (which is
analogous to the N\'eron model of an elliptic curve). In fact, 
the Euler factors of a Drinfeld module at the bad primes are remarkably
similar to those of elliptic curves (ibid.). In
particular, Theorem \ref{cohenmc1} gives the analytic continuation
of the $L$-series of Drinfeld modules and general $A$-modules
\cite{an1}, \cite{boc1} defined over finite extensions of $k$.
\end{rems}

\subsection{The canonical $1$-parameter family of measures associated to
a $\tau$-sheaf}\label{oneparam}

Let $X$, $\underline{\mathcal F}$, etc., be as in Theorem \ref{cohenmc1}.
Let $s=(x,y)\in S_\infty$ and set 
\begin{equation}\label{lseries1}
L(s)=L(x,y):=L(\underline{\mathcal F},s)^{(-1)^{e_X-1}}=\sum_{I}
b_I I^{-s}\,,\end{equation}
where $I$ runs over the ideals of $A$. We know that $L(s)$ has the following
properties:
\begin{enumerate} 
\item The coefficients $\{b_I\}$ belong to $A$.
\item $L(s)$ converges absolutely in some half plane
$\{(x,y)\mid \vert x\vert_\infty \geq t\}$ of $S_\infty$.
\item For each non-negative integer $j$, the power series $L(\pi_\ast^{j}x,-j)$
 is a polynomial in $x^{-1}$ (with coefficients in
${\mathcal O}_\mathbb V$) whose degree is $O(\log(j))$.
\end{enumerate}

\begin{defn}\label{motivic}
Any Dirichlet series $L(s)$ satisfying
the above three properties will be said to be in the {\it motivic class}
$\mathfrak M$
\end{defn}

\noindent
As will be seen by our main result Theorem \ref{main2},
every {\it partial} $L$-series (in the sense
of Definition \ref{partialL} below) associated to a Dirichlet
series in class $\mathfrak M$ will {\em also} be in 
$\mathfrak M$. Note, in particular, 
that Definition \ref{motivic}
definitely does {\em not} require $L(s)$ to have an associated
Euler-product.

\begin{rem}\label{motivated}
We use the adjective ``motivic'' in Definition \ref{motivic} precisely because it is our expectation
that the only general procedure to produce non-trivial Dirichlet
series in class $\mathfrak M$ will be via partial $L$-series of
$\tau$-sheaves.\end{rem}

From now on we let $L(s)=\sum b_I I^{-s}$ be a fixed Dirichlet series
in class $\mathfrak M$.
For $\alpha\in A_{1,\infty}$ let $\delta_\alpha$ be the Dirac measure
concentrated at $\alpha$ as in Definition \ref{dirac}. 

\begin{defn}\label{canonmeasures}
Let $x\in \C_\infty$ have $\vert x\vert_\infty\geq t$ (where
$t$ gives a half-plane of absolute convergence as above). Then we 
define
\begin{equation}\label{canonmeasures2}
\mu_{L,x}:=\sum_I b_Ix^{-\deg_k I} \delta_{\langle I\rangle}\,,
\end{equation}
where $I$ runs over all ideals of $A$. We call $\mu_{L,x}$ {\it the
canonical ($1$-parameter) family of measures} associated to $L(s)$ at $\infty$.
\end{defn}
\noindent
It is clear that, with $x$ chosen as above,
the series for $\mu_{L,x}$ converges to a bounded measure
on $A_{1,\infty}$. As such, its coefficients with respect to any
basis must also be bounded.

Let $v$ be a prime of $A$ and 
let $k_{\sigma,v}$ be as in Equation \ref{ksigmav} with maximal compact
subring $A_{\sigma,v}$. The $v$-adic version of Definition \ref{canonmeasures}
is given next.

\begin{defn}\label{canonmeasures3}
Let $x_v\in \C_v$ with $\vert x_v\vert_v \geq t_v>1$. We define
\begin{equation}\label{canonmesures4}
\mu_{L,x_v}:=\sum_I b_Ix_v^{-\deg I}\delta_{\sigma(I^{s_1})}\,,\end{equation}
where, again, $I$ ranges over {\em all} ideals of $A$.
We call $\mu_{L,x,v}$ the {\it canonical family of measures associated to
$L(s)$ at $v$.}
\end{defn}

\noindent
It is clear that $\mu_{L,x_v}$ also converges to a bounded measure with
the above choice of $x_v$.

The reader should note, of course, that if $v\mid I$, then 
$\sigma(I^{s_1})\not \in A_{\sigma,v}^\ast$.

\section{The main theorem}\label{main}
Let $L(s)=\sum b_I I^{-s}$ continue to be a Dirichlet series in class
$\mathfrak M$.

\subsection{Partial $L$-series}\label{partial}

Let $W=\{w_1,\ldots, w_k\}$ be a finite collection of places of
$\mathbb V$. We explicitly allow at most one element of $W$ to be
the canonical infinite prime of $\mathbb V$ and the
rest are assumed to be finite places. (As in Remark \ref{infiniteplaces},
one may also use {\em all} the infinite primes of $\mathbb V$, and we
leave to the reader the easy modifications necessary to include them.)
Let $E=\{n_1,\ldots,n_k\}$ be a collection of positive integers and let
${\mathfrak w}={\mathfrak w}_{W,E}$ 
be the effective divisor $\sum n_i w_i$ on $C$. Write
${\mathfrak w}={\mathfrak w}_f+{\mathfrak w}_\infty$ where
${\mathfrak w}_f$ consists of the sum over the finite primes and 
${\mathfrak w}_\infty$ is a multiple of $\infty$ (which may be
$0$ if $\infty\not\in W$).

For each finite place $w_i$ of $W$, let ${\mathcal O}_{w_i}={\mathcal O}_{{\mathbb V},w_i}$ be the
associated local ring. If $w_i=\infty$ we {\em define} ${\mathcal O}_\infty=
A_{1,\infty}$. 
Let ${\mathfrak a}=\{\alpha_i\in {\mathcal O}_{w_i}\}$ be a 
collection of elements in these local rings.

\begin{defn}\label{congruence}
Let $I$ be an ideal of $A$. We say that $I\equiv {\mathfrak a}
\pmod {\mathfrak w}$ if $I^{s_1} \equiv \alpha_i\pmod{w_i^{n_i}}$ for
all finite $w_i$ in $\mathfrak w$ and $\langle I \rangle\equiv \alpha_j
\pmod{w_j^{n_j}}$ when $w_j$ is the canonical infinite prime $\infty$.
\end{defn}

Clearly, should the reader desire, one can use the approximation theorem to
replace $\mathfrak a$ by an element of $\mathbb V$. Note also
that if $w_j=\infty$ then we may assume that $\alpha_j$ is a $1$-unit
as otherwise the definition is vacuous.

\begin{defn}\label{partialL}
Let $L(s)=\sum_I b_I I^{-s}$ be a Dirichlet series and let
${\mathfrak a},\,\mathfrak w$ be as above. We set
\begin{equation}\label{partialL2}
L_{{\mathfrak a},\mathfrak w}(s):=\sum_{I\equiv {\mathfrak a}\pmod{\mathfrak w}}
b_II^{-s}\,.\end{equation}
The Dirichlet series
$L_{{\mathfrak a}, \mathfrak w} (s)$ is called the {\it partial Dirichlet
series} associated to $L(s)$ at ${\mathfrak a},{\mathfrak w}$.\end{defn}

The next two statements are then the main results of this paper.

\begin{theorem}\label{main2}
Let $L(s)$ be a Dirichlet series in the motivic
class $\mathfrak M$. Then $L(s)$ analytically continues to
an essentially-algebraic entire function on $S_\infty$. Moreover, any
$v$-adic interpolation of $L(s)$, for $v\in \text{Spec}(A)$,
analytically continues to an essentially-algebraic
entire function on $\C_v^\ast\times S_\sigma$. 
\end{theorem}

\begin{proof}
We first show that $L(s)$ may be analytically continued to essentially
algebraic entire functions at the canonical
place $\infty$ and the finite primes of ${\mathcal O}_\mathbb V$.
Our proof will be to express $L(s)$ as a uniform limit of entire
functions. We then refine this result to establish to establish the
last part of the result.

Let $\mu_{L,x}$ be the canonical family of measures associated to $L(s)$
as given by Equation \ref{canonmeasures2}. Let $t$ be chosen so that
$L(s)$ converges absolutely for $\vert x\vert\geq t$; thus for such
$x$, the series for $\mu_{L,x}$ converges to a bounded measure. Let
$t_1$ be a positive real number less than $t$. The proof proceeds by
analytically extending $L(s)$, $s=(x,y)\in S_\infty$, from the half-plane
of absolute convergence to $\{(x,y)\in S_\infty
\mid \vert x\vert_\infty\geq t_1\}$ in a manner which is uniform in $y$
and $x$. As
$t_1$ is arbitrary the full analytic continuation follows.

Recall that $K_1$ is the extension of $K$ obtained by adjoining
$\langle I \rangle$ for all ideals $I$ of $A$ with maximal compact
subring $A_{1,\infty}$ and maximal ideal
$M_{1,\infty}$. Let $y\in \Zp$ and let $z\in A_{1,\infty}$.
Define 
\begin{equation}\label{zy}
\tilde{z}^y:=\begin{cases} z^y, &\text{for}\, z\equiv 1\pmod{M_{1,\infty}}\\
0,&\text{otherwise\,,}\end{cases}\end{equation}
where $z^y=(1+(z-1))^y$ is computed as before. Clearly the
function $z\mapsto \tilde{z}^y$ is locally-analytic on $A_{1,\infty}$ of
order $1$. Moreover, one checks easily that
$\Vert \tilde{z}^y\Vert_1=1$ for all $y$ so that, as in 
Remark \ref{unifestimates}, we can obtain estimates on its
expansion coefficients which are uniform in $y$.

For $\vert x\vert_\infty\geq t$ we have the basic integral 
\begin{equation}\label{basicintegral}
L(x,y)=\int_{A_{1,\infty}} \tilde{z}^{-y}\, d\mu_{L,x}(z)\,.\end{equation}
The analytic continuation of $L(s)$ proceeds by showing that this integral
extends to $\vert x\vert_\infty \geq t_1$ for all $t_1$ as above.

Let $\{Q_n(z)\}_{n=0}^\infty$ be an orthonormal basis for the
${\mathcal C}(A_{1,\infty},K_{1,\infty})$ consisting of polynomials of degree
$n$, as in Subsection \ref{gentheory}, with measure coefficients
\begin{equation}\label{meascoeffs1}
\left\{b_n=b_n(x)=\int_{A_1} Q_n(z)\,d\mu_{L,x}(z)\right\}
\end{equation}
(N.B., the measure coefficients are now functions of the parameter $x$ in
$\mu_{L,x}$). For $\vert x\vert \geq t$, $b_n(x)$ is bounded as $\mu_{L,x}$ is
bounded.

The next step is to show that $b_n(x)$ is a polynomial
in $x^{-1}$ whose degree is $O(\log (n))$. To see this
write
\begin{equation}\label{qnexp}
Q_n(z)=\sum_{j=0}^n q_{n,j}z^j\,.\end{equation}
Thus
\begin{equation}\label{qnexp2}
b_n(x)=\sum_{j=0}^n q_{n,j}\int_{A_{1,\infty}}z^j\, d\mu_{L,x}(z)\,.
\end{equation}
By definition, $L(x,-j)=\int_{A_{1,\infty}} z^j\,
d\mu_{L,x}(z)$ and, {\em as $L(s)$ is assumed to be motivic}, 
it is a polynomial in $x^{-1}$ (with coefficients in
$K_1$) whose degree is $O(\log (j))$. Therefore $b_n(x)$ is also
a polynomial whose degree (in $x^{-1}$) $d_n$ is $O(\log(n))$. We write
\begin{equation}\label{bni}
b_n(x)=\sum_{i=0}^{d_n} b_{n,i}x^{-i}\,.\end{equation}

In order to apply Theorem \ref{amice1}, we now switch from
using the absolute value $\vert ?\vert_\infty$ to the equivalent
additive functions as in Remark
\ref{addversion}. Let $v_\infty$, $v_{1,\infty}$, $e$, $q_K$ be
as presented there. Then we see that
$L(s)$ converges absolutely for all $\{(x,y)\in S_\infty\mid v_{1,\infty}(x)
\leq -e\log_{q_K}(t)\}$. 
We know that $\mu_{L,x}$ converges to a bounded measure in
this region. Moreover, by choosing $t$ a bit larger, we can assume that
$v_{1,\infty}(\mu_{L,x} (U))\geq 0$ for any compact open
$U$. As the integral of $Q_n(z)$ against such a measure
must also satisfy the {\em same} bounds,
we have
\begin{equation}\label{bni1}
v_{1,\infty}(b_{n,i})+ie\log_{q_K}(t)\geq 0
\end{equation}
uniformly in $n$.
In other words, there is a non-negative constant $C_1$ such that
$$v_{1,\infty}(b_{n,i})\geq - i C_1$$
uniformly in $n$.

Therefore for any $t_1$ sufficiently small we conclude the existence of
a positive constant $C_2$ such that
for $\vert x\vert_\infty\geq t_1$ we have
$v_{1,\infty}(b_n(x))\geq -d_n C_2$.

Now expand $\tilde{z}^y=\sum_n a_{y,n}Q_n(z)$ so that
the basic integral (\ref{basicintegral}) becomes
\begin{equation}\label{basicintegral2}
L(x,y)=\int_{A_{1,\infty}}\tilde{z}^{-y}\,d\mu_{L,x}(z)=\sum_{n=0}^\infty
a_{-y,n}b_n(x)\,.\end{equation}
By Theorem \ref{amice1} and Remark \ref{unifestimates}, we
conclude that $ v_{1,\infty}(a_{-y,n})\geq C_3 n$ uniformly in $y$
for some positive constant
$C_3$. Therefore
$v_{1,\infty}\left(a_{-y,n}b_n(x)\right)\geq C_3n- C_2d_n$.
As $d_n=O(\log (n))$, this goes to $\infty$ with $n$. In particular
the series for the integral converges uniformly in $x$, for
$\vert x\vert_\infty\geq t_1$, and uniformly
in $y$ thus giving the desired analytic continuation.

The analogous $v$-result is even easier since $I^{s_1}$ is always
a $v$-adic integer.
\end{proof}

\begin{theorem}\label{main3} 
Let $L(s)$ be a Dirichlet series in the motivic class $\mathfrak M$.
Then all partial Dirichlet series associated to $L(s)$ are also
in the class $\mathfrak M$.
\end{theorem}

\begin{proof}
We continue with the notations, etc., as in the proof of
Theorem \ref{main2}.
Let ${\mathfrak a}, {\mathfrak w}$ and 
$L_{{\mathfrak a},\mathfrak w}(s)$ be as in Definition \ref{partialL}.
The first two conditions for $L_{{\mathfrak a}, \mathfrak w}(s)$ to be
in the motivic class, as given in Subsection \ref{oneparam}, are easily
checked to follow directly from those of $L(s)$. Thus we only need check the
more subtle third condition.
Let $j$ be a non-negative integer and
\begin{equation}\label{zminusj}
z_{L_{{\mathfrak a}, \mathfrak w}}(x,-j):=L_{{\mathfrak a}, \mathfrak w}(\pi_\ast^jx,-j)\,.\end{equation}
We need to show that $z_{L_{{\mathfrak a}, \mathfrak w}}(x,-j)$ is a polynomial
in $x^{-1}$ whose degree grows like $O(\log (j))$.

To do this we note first that it suffices to handle the case where
$\mathfrak w$ is supported on one place. Indeed, we simply use the
result inductively at each place dividing $\mathfrak w$. 

We begin by writing 
\begin{equation}\label{zminusj2}
z_{L_{{\mathfrak a}, \mathfrak w}}(x,-j)=\sum_{i=0}^\infty c_{j,i} x^{-i}\,.
\end{equation}
By hypothesis we know that $\{c_{j,i}\}\subseteq {\mathcal O}_{\mathbb V}$.
Thus to show $c_{j,i}$ vanishes for some choice of $j$ and
$i$, by Corollary \ref{nonzero2} we need only show
that it has negative degree. We assume first that $\mathfrak w$ is
supported at $\infty$; the $v$-adic case will follow in a similar
fashion. 

Thus suppose that ${\mathfrak w}=n_\infty\infty$ for $n_\infty\geq 1$ and 
${\mathfrak a}=\alpha$ where $\alpha$ is a $1$-unit in $A_{1,\infty}$. 
Let
$\chi$ be the characteristic function of the open subset
$\alpha+M_{1,\infty}^{n_\infty}$ of $A_{1,\infty}$. Obviously,
$\chi$ is locally analytic of order $h=n_\infty$.

We then obtain the following integral representation for the partial
$L$-series $L_{{\mathfrak a},\mathfrak w}(s)$
\begin{equation}\label{basicintegral3}
L_{{\mathfrak a},\mathfrak w}(s)=\int_{A_{1,\infty}}\chi(z) \tilde{z}^{-y}\,
d\mu_{L,x}(z)\,.\end{equation}
Note that the norm of the locally
analytic function of order $n_\infty$, $\chi(z)\tilde{z}^{-y}$, is again 
$1$ and obviously independent of $y$.
 Thus if we write $\chi(z)\tilde{z}^{y}=\sum_n f_{y,n}Q_n(z)$
we find
\begin{equation}\label{basicintegral4}
L_{{\mathfrak a},\mathfrak w}(s)=\sum_n f_{-y,n}b_n(x)\,,\end{equation}
where we have uniform estimates on $v_{1,\infty}(f_{-y,n})$ by
Theorem \ref{amice1}; i.e., there is a positive constant $C_4$, independent
of $y$ (and $-y$), such that
\begin{equation}\label{yn}
v_{1,\infty}(f_{-y,n})\geq C_4 n\,.\end{equation}

Combining Equations \ref{basicintegral4} and \ref{bni} we obtain
the explicit formula
\begin{equation}\label{bniexplicit}
L_{{\mathfrak a}, \mathfrak w}(s)=L_{{\mathfrak a}, \mathfrak w}(x,y)=
\sum_{i=0}^\infty\left(\sum_{n=0}^\infty f_{-y,n}b_{n,i}\right)x^{-i}\,.
\end{equation}
Thus with $j$ as above we obtain 
\begin{equation} \label{bniexplicit2}
c_{j,i}=\pi_\ast^{-ji}\sum_{n=0}^\infty f_{j,n}b_{n,i}\,.
\end{equation}

We now estimate the valuation of 
$\sum_{n=0}^\infty f_{-y,n}b_{n,i}$;
using this estimate in Equation \ref{bniexplicit2} will
allow us to finish the proof.
As the degree $d_n$ of $b_n(x)$ is $O(\log (n))$ we conclude that if
$x^{-i}$ occurs with a non-zero coefficient in $L_{{\mathfrak a},
\mathfrak w}(x,y)$ then $x^{-i}$ must be contributed by those 
$b_n(x)$ where $n$ is {\em exponential} in $i$. Combining this with
with (\ref{yn}) we obtain the {\em fundamental estimate} 
(for $i$ sufficiently large)
\begin{equation}\label{fundestimate}
v_{1,\infty}\left(\sum_{n=0}^\infty f_{-y,n}b_{n,i}\right)\geq C_5 e^{C_6i}\,,
\end{equation}
for positive $C_5$ and $C_6$; this is again independent of $y$.
(The reader should note that $v_{1,\infty}(b_{n,i})$ is greater then
a constant times $i$ and so can be absorbed into the
exponential term as given above.)

Let $b_\mathbb V$ be the number of infinite places of $\mathbb V$.
Using the estimate (\ref{fundestimate}) in
Equation \ref{bniexplicit2}, we see
that the contribution to $\deg_{\mathbb V}(c_{j,i})$
at the place $\infty$ of $\mathbb V$
must be $\leq [{\mathbb V}\colon k]( -C_5e^{C_6i}+ij)$. On
the other hand, as $L_{{\mathfrak a},\mathfrak w}(s)$ has a half-plane of
absolute convergence and coefficients in 
$A$ (as mentioned above),
the discussion in Subsection \ref{domaininfty}
assures us that the other $b_{\mathbb V} -1$ infinite places of
$\mathbb V$ can
contribute at most a positive constant times $ij$ to the
degree. We conclude that for sufficiently large $i$
\begin{equation}\label{fundest}
\deg_{\mathbb V}(c_{j,i})\leq C_7 ij-C_5 e^{C_6i}\,,
\end{equation}
for some positive constant $C_7$.
Elementary estimates show that this expression is then negative for 
$i\gg \log(j)$ which gives the result.

Again the $v$-adic version follows similarly using the fact that the degree
of a principal divisor on a complete curve must vanish. Finally, we
note that the fundamental estimate (\ref{fundestimate}) can also be used
to give another proof of Theorem \ref{main2}.
\end{proof}

\begin{cor}\label{meascompactopen}
Let $L(s)$ be as in {\rm Theorem \ref{main2}}.
Let $U$ be a compact open subset of $A_{1,\infty}$ (resp.\ $A_{\sigma,v}$).
Then $\mu_{L,x}(U)$ (resp.\ $\mu_{L,{x_v}}(U)$) is a polynomial in
$x^{-1}$ (resp.\ $x_v^{-1}$).
\end{cor}

\begin{proof}
We can suppose $U$ is of the form $\alpha + M_{1,\infty}^j$ (resp.\
$\alpha+ M_{\sigma,v}^j$).
Then $\mu_{L,x}(U)$ (resp.\ $\mu_{L,x_v}(U)$) is the value of the
associated partial $L$-series at $y=0$ (resp.\ $y_\sigma=0$). Thus
the result follows immediately from Theorem \ref{main2}.\end{proof}

The above proof of the analytic continuation depends crucially
on the fact that $z \mapsto \tilde{z}^y$ is locally analytic. If one had
{\em any} other type of locally analytic endomorphism of the
group of $1$-units, then the proof would {\em automatically} work for it also.
Therefore the following question is of great interest.

\begin{question}\label{fundques}
Let $K=\Fr((1/T))$ and let $U_1$ be the group of $1$-units of $K$. 
Does there exist a locally-analytic endomorphism $f\colon U_1\to U_1$ 
which is not of the form $u\mapsto u^y$ for some $p$-adic integer $y$?
\end{question}

It is reasonable to expect a negative answer to Question \ref{fundques}
but we certainly have no proof of this as of this writing.

\section{Complements}\label{complements}
Let $k$ be as before and let $\mathcal E$ be an $A$-module \cite{an1}
\cite{boc1} defined over a finite extension $k_1$ of $k$. For each
finite prime $v$ of $A$ we can define the $v$-adic Tate module
of $\mathcal E$. Using the invariants of inertia and the Frobenius
element at another finite prime $\mathfrak p$ not dividing
$v$ one obtains the local Euler-factors of the $L$-series
of $E$ at $\mathfrak p$ (as mentioned above). 
One can also use the same definitions when $\mathfrak p$ is an {\em infinite}
place. Remarkably, as in \cite{ga1}, one still obtains polynomials
with coefficients in $A$ and which are independent of $v$. For instance,
in the case that $\mathcal E$ is a Drinfeld module (or just
uniformizable) one finds that the coefficients are actually in
$\Fr$ by using the associated lattice (in fact, by Gardeyn \cite{ga2}
this essentially characterizes uniformizable $A$-modules).

As mentioned in \cite{boc1} (also \cite{go6}) these factors at $\infty$ should give
rise to trivial-zeroes for the special polynomials of $L({\mathcal E},s)$,
 and thus
also $L({\mathcal E},s)$ 
itself. In fact, one expects this to ultimately be a completely
general phenomenon for all $\tau$-sheaves.
The trivial zeroes for the $v$-adic functions will be given by
the Euler factors in the special polynomials lying over $v$ (which are removed when one interpolates
$v$-adically).
In the case where $\mathcal E$ is uniformizable
the trivial zeroes resemble the classical case. However, when $\mathcal E$ is 
arbitrary they represent something new.
As of this writing these trivial zeroes have an 
extremely mysterious effect on the collection of all
zeroes. Indeed, because of the 
non-Archimedean nature of
$S_\infty$ and the uniform continuity of $L(s)$, trivial
zeroes influence nearby zeroes (called ``near-trivial zeroes'') as
in \cite{go6}. One therefore wants to find the set of zeroes (called
``critical zeroes'') which are not so influenced. This seems very
hard at the moment. Classical theory suggests that it may require
deeper understanding of the connections with modular forms as in
\cite{boc2} (see also \cite{go7}).

\begin{rem}\label{dinesh}
We finish by explaining how elementary estimates can be used to imply
logarithmic growth in general in the rank $1$ case (as is evident
in the explicit calculation
given in \cite{th1}).  Let $\mathcal K$ be any field of characteristic $p$ 
and let
$W\subseteq \mathcal K$ be a finite additive subgroup of order $p^{m_W}$. Define
$s_i(W):=\sum_{w\in W}w^i$ for non-negative integers $i$. Then
one has
\begin{equation}\label{dinesh2}
0\leq i< p^{m_W}-1\Rightarrow s_i(W)=0\,.\end{equation}
The proof, \`a la Carlitz,
goes as follows. Put $e_W(z):=\prod_{w\in W}(z-w)$ which, by
standard arguments, is an additive polynomial with constant derivative
$\lambda_W\neq 0$. Logarithmic differentiation then implies
$$\lambda_W/e_W(z)=\sum_{w\in W}(z-w)^{-1}=\sum_{w\in W}\frac{1}{z}(1-w/z)^{-1}\,.$$ 
Using the geometric series,
one finds that the coefficient of $z^{-j}$ is $s_{j-1}(W)$. On the
other hand, $\lambda_W/e_W(z)$ has a zero of order $p^{m_W}$ at $\infty$
and the result follows. Applying this to the sums over finite subgroups which
arise in the rank $1$ case immediately gives the logarithmic growth
of special polynomials.

Note that Equation \ref{dinesh2} includes the very well-known
result $\dis \sum_{a\in \Fr}a^i=0$ for $0\leq i<r-1$.
\end{rem}


\begin{thebibliography}{3333334}
\bibitem[Am1]{am1}{\sc Y.\ Amice:} Interpolation $p$-adique, {\it Bull.\
Soc.\ Math.\ France} {\bf 92} (1964) 117-180.

\bibitem[Am2]{am2} {\sc Y.\ Amice:} Duals, in: {\it Conference
on $P$-adic Analysis, (1978: Nijmegen, Netherlands}),
Roomsch-Katholieke Universiteit (Nijmegen, Netherlands),
Report-Mathematische instituut, Katholieke universiteit; 7806.

\bibitem[An1]{an1}  {\sc G.\  Anderson:}  $t$-motives, 
{\it Duke Math. J.} {\bf 53} (1986) 457-502.

\bibitem[Boc1]{boc1} {\sc G.\ B\"ockle:} 
Global $L$-functions over function fields,
{\it Math.\ Ann.} {\bf 323} (2002) 737-795.

\bibitem[Boc2]{boc2} {\sc G.\ B\"ockle:} An Eichler-Shimura isomorphism over 
function fields between Drinfeld modular forms and cohomology classes of
crystals,(preprint, available at\\ \verb+http://www.math.ethz.ch/~boeckle/+).

\bibitem[BP1]{bp1} {\sc G.\ B\"ockle, R.\ Pink:} A cohomological 
theory of crystals over function fields, (in preparation).

\bibitem[Ca1]{ca1} {\sc L.\ Carlitz:} On certain functions connected with 
polynomials in a Galois field, {\it Duke Math.\ J.} {\bf 1} (1935) 137-168.

\bibitem[Co1]{co1} {\sc K.\ Conrad}: The digit principle, {J.\ Number Theory}
{\bf 84} (2000) 230-257.

\bibitem[DV1]{dv1} {\sc J.\ Diaz-Vargas}: Riemann 
hypothesis for $\Fp[T]$, {\it J.\ Number Theory} {\bf 59} (1996) 313-318

\bibitem[Dr1] {dr1} {\sc V.G.\ Drinfeld:}  
Elliptic modules, Math.\ Sbornik {\bf 94} 
(1974) 594-627, English transl.:  {\it Math.\ U.S.S.R.\ Sbornik} {\bf 23} 
(1976) 561-592.

\bibitem[Ga1]{ga1} {\sc F.\ Gardeyn:} A Galois criterion for good
reduction of $\tau$-sheaves, {\it J.\ Number Theory} {\bf 97} (2002) 447-471.

\bibitem[Ga2]{ga2} {\sc F.\ Gardeyn:} New criterion for Anderson 
uniformizability
of $T$-motives, (preprint, available at 
\verb+ http://www.math.ethz.ch/~fgardeyn/FG/index.shtml+).

\bibitem[Ga3]{ga3} {\sc F.\ Gardeyn:} The structure of 
analytic $\tau$-sheaves, {\it J.\ Number Theory} {\bf 100} (2003) 332-362.

\bibitem[Go1]{go1} {\sc D.\ Goss}: $\pi$-adic Eisenstein Series 
for Function Fields, {\it Compositio Math.} {\bf 41} (1980) 3-38.

\bibitem[Go2]{go2} {\sc D.\ Goss}: Modular Forms for ${\bf F}_r[T]$, 
{\it J.\ Reine Angew.\ Math.} {\bf 317} (1980) 16-39.

\bibitem[Go3]{go3} {\sc D.\ Goss}: Some integrals associated to modular
forms in the theory of function fields, in: 
{\it The Arithmetic of Function Fields} 
(eds:  D. Goss et al) de Gruyter (1992) 227-251.

\bibitem[Go4]{go4} {\sc D.\ Goss}:
{\it Basic Structures of Function Field Arithmetic},
Springer-Verlag, Berlin (1996).

\bibitem[Go5]{go5} {\sc D.\ Goss}: A Riemann hypothesis for 
Characteristic $p$ $L$-functions, {\it J.\ Number Theory} {\bf 82}
(2000) 299-322.

\bibitem[Go6]{go6} {\sc D.\ Goss}: The impact of the infinite primes on the
Riemann hypothesis for characteristic $p$ valued $L$-series, 
in: {\it Algebra, Arithmetic and Geometry with Applications}
(Eds: Christensen et al) Springer (2004) 357-380.

\bibitem[Go7]{go7} {\sc D.\ Goss}: 
Can a Drinfeld module be modular? J.\ Ramanujan Math.\ Soc.\
{\bf 17} No.\ 4 (2002) 221-260. 

\bibitem[Go8]{go8} {\sc D.\ Goss}: 
Fourier Series, Measures and Divided Power Series in the Theory of
Function Fields, $K$-Theory {\bf 1} (1989) 533-555.

\bibitem[Hay1]{hay1} {\sc D.\ Hayes:} A brief introduction 
to Drinfeld modules, 
in: {\it The Arithmetic of Function Fields} (eds. D. Goss et al) de Gruyter
(1992) 1-32.

\bibitem[SJ1]{sj1} {\sc S.\ Jeong:} A comparison of the Carlitz and 
digit derivative bases in function field arithmetic, {\it J.\ Number
Theory} {\bf 84} (2000) 258-275.

\bibitem[SJ2]{sj2} {\sc S.\ Jeong:} Continuous linear endomorphisms 
and difference equations over the completion of $\Fq[T]$, {\it J.\
Number Theory} {\bf 84} (2000) 276-291.

\bibitem[SJ3]{sj3} {\sc S.\ Jeong:} Hyperdifferential operators and 
continuous functions on function fields, {\it J.\ Number Theory} {\bf 89}
(2001) 165-178.

\bibitem[SJ4]{sj4} {\sc S.\ Jeong:} Digit derivatives and 
application to zeta measures, (preprint).

\bibitem[Ka1]{ka1} {\sc I.\ Kaplansky}: The Weierstrass theorem in fields
with valuations, {\it Proc.\ Amer.\ Math.\ Soc.} {\bf 1} (1950) 356-357.

\bibitem[Ko1]{ko1} {\sc N.\ Koblitz}: {\it $p$-adic Analysis: a Short
Course on Recent Work}, London Math.\ Soc.\ Lecture Note Series {\bf 46},
Cambridge Univ.\ Press (1980)

\bibitem[R1]{r1} {\sc B.\ Riemann:} Ueber die Anzahl der Primzahlen
unter einer gegebenen Gr\"osse, {\it Monatsberichte der Berliner Akademie}
(1859); {\it Gesammelte Werke,} Teubner, Leipzig
(1892).

\bibitem[Ros1]{ros1} {\sc M.\ Rosen}: {\it Number Theory in
Function Fields}, Springer 2002.

\bibitem[Sh1]{sh1} {\sc J.\ Sheats}: The Riemann hypothesis for the Goss
zeta function for $\Fq[T]$, {\it J.\ Number Theory} {\bf 71} (1998) 121-157.

\bibitem[TW1]{tw1} {\sc Y.\ Taguchi, D.\ Wan:} 
$L$-functions of $\varphi$-sheaves and
Drinfeld modules, {\it J. Amer. Math. Soc.} {\bf 9} (1996) 755-781.

\bibitem[TW2]{tw2} {\sc Y.\ Taguchi, D.\ Wan:} Entireness of 
$L$-functions of $\varphi$-sheaves on affine complete 
intersections, {\it J.\ Number Theory} {\bf 63} (1997) 170-179.

\bibitem[Th1]{th1} {\sc D.\ Thakur:} Zeta measure associated to 
$\mathbb F_q[T]$, {\it J. Number Theory} {\bf 35} (1990) 1-17.

\bibitem[vR1]{vr1} {\sc A.\ C.\ van Rooij:} {\it Non-Archimedean Functional
Analysis}, Marcel Dekker (1978).

\bibitem[Wa1]{wa1} {\sc D.\ Wan:} On the Riemann hypothesis for the
characteristic $p$ zeta function,  {\it J.\ Number Theory} {\bf 58} (1996)
 196-212.

\bibitem[Ya1]{ya1} {\sc Z.\ Yang:} Locally analytic functions over
completions of $\Fr[U]$, {\it J.\ Number Theory} {\bf 73} (1998) 451-458.

\bibitem[Ya2]{ya2} {\sc Z.\ Yang:} A note on zeta measures over
function fields, {\it J.\ Number Theory} {\bf 90} (2001) 89-112.

\end{thebibliography}
\end{document}